\documentclass[12pt]{article}
\usepackage[latin1]{inputenc}
\usepackage{amsmath}
\usepackage{amsfonts}
\usepackage{amssymb}
\usepackage{color}
\usepackage{graphicx}
\newtheorem{theorem}{Theorem}[section]

\newcommand{\qed}{\hspace{\stretch{3}}$\square$\\[1.8ex]}

\newtheorem{lemma}[theorem]{Lemma}
\newtheorem{cor}[theorem]{Corollary}

\newtheorem{teo}[theorem]{Theorem}
\newtheorem{remark}[theorem]{Remark}

\def\R{\mathbb{R}}

\def\R{\mathbb{R}}

\def\M{\mathcal{M}}

\def\F{\mathcal{F}}

\def\l{\langle}

\def\r{\rangle}

\newcommand{\proof}{{\bf Proof: }}

\title{\bf  Supercritical super-Brownian
motion with a general branching mechanism and travelling waves}
\author{{\sc A.E. Kyprianou\footnote{Department of Mathematical Sciences, University of Bath,
Claverton Down, Bath, BA2 7AY, U.K.} \ \  R.-L.
Liu\footnote{Department of Mathematics,
Nanjing University, Nanjing 210093,
 P. R. China.}  \ \
A. Murillo-Salas\footnote{Departamento de Matem\'aticas, Universidad de Guanajuato,
Jalisco s/n, Mineral de Valenciana,
Guanajuato, Gto. C.P. 36240, M\'exico.} \ \
  Y.-X.
Ren\footnote{LMAM School of Mathematical Sciences,
Peking University, Beijing 100871,
 P. R. China.} 
}
}

\begin{document}

\maketitle
\begin{abstract}
We consider the classical problem of existence, uniqueness and asymptotics of monotone solutions to the travelling wave equation associated to the
parabolic semi-group equation of a super-Brownian motion with a general branching mechanism. Whilst we are strongly guided by the probabilistic
reasoning of Kyprianou \cite{Kyprianou2004} for branching Brownian motion, the current paper offers a number of new insights. Our analysis incorporates
the role of Seneta-Heyde norming which, in the current setting, draws on classical work of Grey \cite{Grey}. We give a {\it pathwise} explanation of
Evans' immortal particle picture (the spine decomposition) which uses the Dynkin-Kuznetsov $\mathbb{N}$-measure as a key ingredient. Moreover, in the
 spirit of Neveu's stopping lines we make repeated use of Dynkin's exit measures. Additional complications arise from the general nature of the branching
 mechanism. As a consequence of the analysis we also offer an exact $X(\log X)^2$ moment dichotomy for the almost sure convergence of the so-called derivative martingale at its critical parameter to a non-trivial limit. This differs to the case
 of branching Brownian motion, \cite{Kyprianou2004}, and branching random walk, \cite{BK}, where a moment `gap' appears in the necessary and sufficient
conditions.

\bigskip


\bigskip


\end{abstract}

\section{Introduction}

Suppose that $X=\{X_t: t\geq 0\}$ is a (one-dimensional) $\psi$-super-Brownian motion with general branching mechanism $\psi$ taking the form
\begin{equation}
\psi(\lambda) = -\alpha \lambda + \beta\lambda^2 + \int_{(0,\infty )} (e^{-\lambda x} - 1 + \lambda x )\nu({\rm d}x),
\label{branch-mech}
\end{equation}
for $\lambda \geq 0$ where $\alpha = - \psi'(0^+)\in(0,\infty)$, $\beta\geq 0$ and
$\nu$ is a measure concentrated on $(0,\infty)$ which satisfies
$\int_{(0,\infty)}(x\wedge x^2)\nu({\rm d}x)<\infty$. 
Let
$\mathcal{M}_F(\mathbb{R})$ be  the space of finite measures on $\mathbb{R}$
and note that $X$ is a $\mathcal{M}_F(\mathbb{R})$-valued Markov
process under $\mathbb{P}_\mu$ for each $\mu\in\mathcal{M}_F(\mathbb R)$, where $\mathbb{P}_\mu$ is law of $X$ with initial configuration $\mu$. We shall use standard inner product notation, for
$f\in C_b^+(\mathbb{R})$ and $\mu\in\mathcal{M}_F(\mathbb{R})$,
\[
\langle f , \mu\rangle =  \int_{\mathbb{R} }f(x)\mu({\rm d}x).
\]
Accordingly we shall write $||\mu|| = \langle 1,\mu \rangle$.

The existence of our class of superprocesses is guaranteed by \cite{D, Dyn1993, Dyn2002}. The following  standard result from the theory of superprocesses  describes the evolution of $X$ as a Markov process.
For all $f\in C_b^+(\mathbb{R})$, the space of positive, uniformly bounded, continuous functions on $\mathbb{R}$,   and $\mu\in\mathcal{M}_F(\mathbb{R})$,
\begin{equation}
 -\log\mathbb{E}_\mu(e^{- \langle f, X_t\rangle}) =  \int_{\mathbb{R} }u_f(x, t)\mu({\rm d}x), \, \mu\in\mathcal{M}_F(\mathbb{R}), \, t\geq 0,
 \label{prePDE}
\end{equation}
where  $u_f(x,t)$ is the unique positive solution to the evolution equation for $x\in\mathbb{R}$ and $t>0$
      \begin{equation}\label{PDE}
       \dfrac{\partial }{\partial
       t}u_f(x,t)=\dfrac{1}{2}\dfrac{\partial^2 }{\partial
       x^2}u_f(x,t)-\psi(u_f(x,t)),
       \end{equation}
with initial condition $u_f(x,0) = f(x)$. The reader is referred to Theorem 1.1 of  Dynkin \cite{Dyn1991}, Proposition 2.3 of Fitzsimmons \cite{Fitz} and Proposition 2.2 of Watanabe \cite{watanabe1968} for further details; see also Dynkin  \cite{Dyn1993, Dyn2002} for a general overview.
The analogous object to (\ref{PDE}) for branching Brownian motion is called the Fisher-Kolmogorov-Petrovski-Piscounov (FKPP) equation and hence in the current setting we name (\ref{PDE}) the FKPP equation for $\psi$-super-Brownian motion.

Recall that the total mass of the process $X$ is a continuous-state branching process with branching mechanism $\psi$. Since there is no interaction between spatial motion and branching
we can  characterise the $\psi$-super-Brownian into the
categories of supercritical, critical and subcritical accordingly
with the same categories for continuous-state branching processes.
Respectively, these cases correspond to $\psi'(0^+)<0$, $\psi'(0^+)=0$
and $\psi'(0^+)>0$.  The class of $\psi$-super-Brownian motions described above are necessarily supercritical. Such processes may exhibit explosive behaviour, however, under the conditions assumed above, $X$ remains finite at all positive times.
We insist moreover that $\psi(\infty)=\infty$ which means that with positive probability the event $\lim_{t\uparrow\infty}||X_t||=0$ will occur. Equivalently this means that the total mass process does not have monotone increasing paths; see for example the summary in Chapter 10 of Kyprianou \cite{K}. 
The probability of the
event
\[
\mathcal{E}:= \{\lim_{t\uparrow \infty} ||X_t|| =0\}
\]
 is described in terms of the largest root, say $\lambda^*$, of
the equation $\psi(\lambda)=0$. Note that it is known (cf. Chapter 8 of
\cite{K}) that $\psi$ is strictly convex with $\psi(0)=0$ and
hence since $\psi(\infty)=\infty$ and $\psi'(0^+)<0$ it follows that
there are exactly two roots in $[  0,\infty)$, one of which is always $0$. For $\mu\in
\mathcal{M}_F(\mathbb{R})$ we have
\begin{equation}
 \mathbb{P}_\mu (\lim_{t\uparrow \infty} ||X_t|| =0) = e^{-\lambda^* ||\mu||}.
 \label{extinguishpr}
\end{equation}
In this article we shall also assume on occasion that
  \begin{equation}
 \int^\infty\frac{1}{\sqrt{\int_{\lambda^*}^\xi \psi(u) {\rm d}u}}{\rm d}\xi<\infty. 
 \label{A3}
 \end{equation}
 This condition is equivalent to requiring that $\int^\infty \left(\int_{0}^\xi ( \psi(u)+\alpha u) {\rm d}u\right)^{-1/2}{\rm d}\xi<\infty$.
 In combination with additional assumptions on $\psi$ given above, (\ref{A3}) has a number of implications for the underlying superprocess. Firstly, if we denote  by $\mathcal{R}$ the smallest closed set in $\mathbb{R}$
 such that ${\rm supp} X_t\subseteq\mathcal{R}$ for all $t\geq 0$, then Sheu \cite{Sheu1997} shows that
for all $\mu\in\mathcal{M}_F(\mathbb{R})$ with compact support,
 \[
 \mathbb{P}_\mu(\mathcal{R}\text{ is compact})=e^{-\lambda^* ||\mu||}.
 \]
Secondly (\ref{A3}) implies that $\int^\infty 1/\psi(\xi){\rm d}\xi<\infty$ (cf. \cite{Sheu1997}) which in turn guarantees that the event
$
\mathcal{E}
$
 agrees with the event of {\it extinction}, namely $\{\zeta<\infty\}$ where
\[
 \zeta = \inf\{t> 0 : ||X_t || = 0\}.
\]
 Note that (\ref{A3}) cannot be satisfied for branching mechanisms which belong to bounded variation spectrally positive L\'evy processes.

\bigskip

Our primary concern in this paper will be to look at monotone travelling wave solutions to the FKPP equation (\ref{PDE}). Specifically, we are interested in non-increasing solutions to (\ref{PDE}) of the form $\Phi_c(x - ct)$, where $\Phi_c\geq 0$ and  $c\geq 0$ is the wave speed. That is to say $\Phi_c$ solves
\begin{equation}
\frac{1}{2}\Phi_c'' + c\Phi'_c - \psi(\Phi_c) = 0.
\label{travwave}
\end{equation}
Moreover, for technical reasons which will  become clear later, we shall be interested in the case that
\[
\Phi_c(-\infty) = \lambda^* \mbox{ and }\Phi_c(+\infty ) = 0.
\]
Henceforth we shall say that any solution to (\ref{travwave}) which respects the aforementioned conditions of non-negativity, monotonicity and connecting the points $\lambda^*$ at $-\infty$ to $0$ at $+\infty$ is a {\it travelling wave with wave speed $c$}.

The case that $\psi(\lambda) = -\alpha\lambda + \beta\lambda^2$ corresponds to quadratic branching which has received a great deal of attention in the past in connection with branching Brownian motion. Indeed, note that $\lambda^* = \alpha/\beta$  and hence a simple calculation shows that $\phi_c(x): = 1-(\beta/\alpha) \Phi_c (x)$ satisfies
\begin{equation}
\frac{1}{2}\phi_c'' + c\phi'_c +\alpha (\phi_c^2 - \phi_c)
 =0,
 \label{classicalKPP}
\end{equation}
with $\phi_c(-\infty) = 0$ and $\phi_c(\infty) = 1$.
Starting with Kolmogorov et al. \cite{Kol} and Fisher \cite{Fish} there exists a variety of analytical treatments of
travelling wave equations similar to (\ref{travwave}) and (\ref{classicalKPP}). We name but a few, for example
Aronson and Weinberger \cite{AW}, Fife and McLeod \cite{FM}, Bramson \cite{Bra}, Lau \cite{Lau}, Pinsky \cite{P},
Kametaka \cite{Kam}, Volpert et al. \cite{Volpert}. Our interest however lies more in the probabilistic direction.
There exists a suite of literature which gives a probabilistic handling of (\ref{classicalKPP}); see  \cite{McKean, neveu, Bra, chauvin, DN,  SHarris1999, Kyprianou2004}. Key to all of these papers is the relationship between the travelling wave equation and two types of martingales commonly referred to as {\it additive} and {\it multiplicative} martingales. Our objective in this paper is to show that many of the known probabilistic ideas can be adapted, subject to the use of appropriate alternative technologies,   to handle  (\ref{travwave}). In particular we shall largely work with Dynkin exit measures as well a new {\it pathwise} version of Evans' immortal particle decomposition of our $\psi$-super-Brownian motion.

The remainder of the paper is structured as follows. In the next section we state our main results. These pertain to a complete existence, uniqueness and asymptotics result for the travelling wave equation. In special cases, it is possible to give more explicit details concerning the form of the solutions to the travelling wave equation in terms of martingale limits. For this reason, part of our main results includes some martingale convergence theorems. One of our martingale results, concerning the question of convergence to a non-trivial limit of the so-called derivative martingale, offers a moment dichotomy which has not been previously achieved for the analogous martingales in the case of branching Brownian motion and branching random walks. In Section \ref{BEM} we examine certain Dynkin exit measures which   will be key to later analysis. The remaining sections are dedicated to the proofs of the main results with the exception of Section \ref{spinesection} which provides
  the new pathwise spine (or immortal particle) decomposition that features heavily in the proofs.

On a final note, we mention that the condition (\ref{A3}) appears to be a natural sufficient condition under which to perform all of our analysis. This will become clear later on through several of the preparatory results. We refrain from imposing this condition throughout the paper however (in favour of stating it when required) as a number of the mathematical tools we appeal to, which are of intrinsic interest on their own, still have meaning without it.

\section{Main results}
Our first result gives us a very general characterisation of the existence, uniqueness and asymptotics of non-negative travelling waves solving (\ref{travwave}). Subsequently we give moment conditions under which some of the quantities involved can be explicitly identified. For convenience we write $\underline{\lambda}=\sqrt{-2\psi'(0^+)}$ and for each $\lambda\in\R$ define
\begin{equation}
c_\lambda=-\psi^\prime(0^+)/\lambda+\lambda/2.
\label{clambda}
\end{equation}
 Note that for $\lambda\in(0, \underline{\lambda}]$, $c_\lambda$ has range $[\underline{\lambda},\infty)$. In particular $c_{\underline{\lambda}} = \underline{\lambda}$.  We shall also write $\mathbb{P}$ as shorthand for $\mathbb{P}_{\delta_0}$ with corresponding expectation operator given by $\mathbb{E}$.

\begin{theorem}
\label{main-1}
$\mbox{ }$
\begin{itemize}
\item[(i)] If (\ref{A3}) holds then  no travelling waves exist with wave speed $c$ if
$c\in[ 0,c_{\underline{\lambda}})$.
\item[(ii)] A  travelling wave exists with wave speed  $c$ if $c\geq c_{\underline{\lambda}}$.
In particular for
$\lambda\in(0,\underline{\lambda}]$ there exists a
travelling wave with wave speed $c_\lambda$ which may be written in
the form
         \begin{equation}\label{subcriticasolution}
         \Phi_{c_\lambda}(x)=-\log \mathbb{E}\left[e^{-e^{-\lambda x}\Delta(\lambda)}\right],
         \end{equation}
   where $\Delta(\lambda)$ is a non-negative random variable such that $\{\Delta(\lambda)=0\}$ agrees with $\mathcal{E}$,  $\mathbb{P}$-almost surely.

\item[(iii)] Suppose that $\lambda\in (0,\underline{\lambda}]$.  Then, up to an additive shift in its argument, there is a unique travelling wave at speed $c_\lambda$.
    \item[(iv)]   Moreover, when
    $\lambda\in(0,\underline{\lambda}]$, there exists some constant $k_\lambda\in(0,\infty)$ and
    a slowly varying function $L_\lambda: (0,\infty)\mapsto (0,\infty)$
    such that
         \begin{equation}\label{subcriticalasym}
         \lim_{x\rightarrow\infty}  \frac{\Phi_{c_\lambda}(x)}{e^{-\lambda x} L_\lambda(e^{-\lambda x})}=k_\lambda.
         \end{equation}
\end{itemize}
\end{theorem}

The conclusions given in the above theorem conform largely to what is understood for the classical FKPP equation  and the folk law of extensions thereof (cf. \cite{Kol, Fish, AW, FM, Lau, Bra, P, Kam, Volpert}) in the sense of existence, uniqueness and asymptotic decay. One might note however that the general form of the slowly varying correction to the exponential decay given in  (\ref{subcriticalasym}) may well be a new result that is not to be found in the PDE literature. There are instances however where, under further assumptions the function $L_\lambda (z)$ is known to behave as $-\log z$ as $z\downarrow 0$. This conclusion will also appear shortly in the forthcoming Theorem \ref{main-2}.

As alluded to above, the next two main theorems make a clearer
statement about the quantities $\Delta(\lambda)$ and $L_\lambda$
when we impose additional assumptions. To do this, we need to
introduce two families  of $\mathbb{P}$-martingales with respect to
the natural filtration $\F_t:=\sigma(X_u; u\leq t)$. The first such family of  martingales is
identified in the following Lemma.
\begin{lemma}
 The process  $W(\lambda) = \{W_t(\lambda) : t\geq 0\}$ where $\lambda\in\mathbb{R}$ and
              \begin{equation}\label{mtg}
              W_t(\lambda):=e^{-\lambda c_\lambda t}\langle e^{-\lambda\cdot},X_t\rangle,\,\,t\geq0,
              \end{equation}
  is a martingale.
  \end{lemma}

\proof
The proof appeals to a classical technique which we briefly outline. Define for each $x\in\mathbb{R}$, $g\in C^+_b(\mathbb{R})$ and $\theta, t\geq 0$, $u_g^\theta(x,t) = - \log \mathbb{E}_{\delta_x}(e^{-\theta\langle g, X_t\rangle})$ and note that, with limits understood as $\theta\downarrow 0$,  $u_g(x,t)|_{\theta =0} =0$ and  $v_g(x, t): =\mathbb{E}_{\delta_x}(\langle g, X_t\rangle) = \partial u_g^\theta(x,t)/\partial\theta |_{\theta = 0}$. Differentiating in $\theta$ in (\ref{PDE}) shows that $v_g$ solves the equation
   \begin{equation}\label{PDE2}
       \dfrac{\partial }{\partial
       t}v_g(x,t)=\dfrac{1}{2}\dfrac{\partial^2 }{\partial
       x^2}v_g(x,t)-\psi'(0^+)v_g(x,t),
       \end{equation}
       with $v_g(x,0) = g(x)$.  Note that classical Feynman-Kac theory tells us that (\ref{PDE2}) has a unique solution and it is necessarily equal to $\Pi_x(e^{-\psi'(0^+) t} g(\xi_t))$ where $\{\xi_t : t\geq 0\}$ is a Brownian motion issued from $x\in\mathbb{R}$  under the measure $\Pi_x$.
The above procedure also works for $g(x) =  e^{-\lambda x}$ in which
case we easily conclude that for all $x\in\mathbb{R}$ and $t\geq 0$,
$e^{-\lambda c_\lambda t}\mathbb{E}_{\delta_x}(\langle e^{-\lambda
\cdot}, X_t \rangle) = e^{-\lambda x}$. Finally, the martingale
property follows using the previous equality together with the
Markov branching property associated with $X$. \qed

\bigskip

Note that $W(\lambda)$ is a nonnegative martingale and
therefore converges almost surely.  As a corollary to the above
lemma, we may describe the second family of martingales we are
interested  in by taking the negative derivative in $\lambda$ of
$W(\lambda)$. Note that this produces a signed martingale which does
not necessarily converge almost surely.

\begin{cor}
The process
$\partial W(\lambda): = \{\partial W_t(\lambda),t\geq0\}$, where $\lambda\in\mathbb{R}$ and
              \begin{equation}\label{derivativemartingale}
              \partial W_t(\lambda):=-\frac{\partial}{\partial\lambda}W_t(\lambda)
                                    =\l(\lambda t+\cdot)e^{-\lambda(c_\lambda t+\cdot)},X_t\r,\,\,t\geq0,
              \end{equation}
is also a martingale.
\end{cor}

It turns out that the convergence of both these martingales in the appropriate sense is important to give a more precise characterization of the limit $\Delta(\lambda)$ and the normalizing sequence $L_\lambda$ in Theorem \ref{main-1}. The following theorem contains the necessary information.

\begin{theorem}\label{mgcgce}$\mbox{ }$
\begin{itemize}
\item[(i)] The almost sure limit of $W(\lambda)$, denoted by $W_\infty(\lambda)$,
is also an $L^1(\mathbb{P})$-limit if and only if
$|\lambda|<\underline{\lambda}$ and
 $$\int_{[1,\infty)} r(\log r) \nu({\rm d}r)<\infty.$$ When $W_\infty(\lambda)$ is an $L^1(\mathbb{P})$-limit the event $\{W_\infty(\lambda)>0\}$ agrees with $\mathcal{E}^c$, $\mathbb{P}$-almost surely. Otherwise, when it is not an $L^1(\mathbb{P})$-limit, its limit is identically zero.

\item[(ii)] Assume that (\ref{A3}) holds.
The martingale $\partial W(\lambda)$ has an almost sure non-negative limit
when $|\lambda|\geq \underline\lambda$ which is identically zero
when $|\lambda|>\underline\lambda$ and when $|\lambda| =
\underline{\lambda}$ its limit is almost surely
strictly positive on  $\mathcal{E}^c$ if and only if
$$\int_{[1,\infty)} r (\log r)^2\nu({\rm d}r)<\infty.$$
\end{itemize}
\end{theorem}

\begin{remark}\rm Note that other similar theorems exist for derivative martingales in the branching random walk, \cite{BK}, and branching Brownian motion, \cite{Kyprianou2004}. In those cases however, an exact dichotomy for convergence to a non-zero limit in the critical regime was not achieved unlike the case here.
\end{remark}

\bigskip

We may now turn to our final main theorem which is a refinement of Theorem \ref{main-1} under additional assumptions. For convenience we write $W_\infty(\lambda)$ and $\partial W_\infty(\lambda)$ for the martingale limits (when it exists in the latter case).
\begin{teo} \label{main-2}
Assume  (\ref{A3}).
\begin{itemize}
\item[(i)] Suppose that $\int_{[1,\infty)} r(\log r)\nu({\rm d}r)<\infty$ and
   $\lambda\in(0,\underline{\lambda})$.
   Then, up to an additive constant in its argument, the
   travelling wave solution $\Phi_{c_\lambda}$ to (\ref{travwave}) is
   given by
         \begin{equation}\label{subcriticasolution}
         \Phi_{c_\lambda}(x)=-\log \mathbb{E}\left[e^{-e^{-\lambda x}W_\infty(\lambda)}\right],
         \end{equation}
   and $L_\lambda(x)\sim 1$ as $x\downarrow 0$. 

\item[(ii)] Suppose that  $\int_{[1,\infty)} r(\log r)^2\nu({\rm d}r)<\infty$
and
   $\lambda=\underline{\lambda}$.   Then, the critical
   travelling wave solution $ \Phi_{\underline\lambda}$ to (\ref{travwave}) is
   given by
         \begin{equation}\label{subcriticasolution}
         \Phi_{\underline\lambda}(x)=-\log \mathbb{E}\left[e^{-e^{-\underline\lambda x}\partial W_\infty(\underline\lambda)}\right].
         \end{equation}
  Moreover, $L_{\underline{\lambda}}(x)\sim -\log x$ as $x\downarrow 0$.
   \end{itemize}
   \end{teo}
\begin{remark}\rm
Note that $$\int^1_0r^{-2}\psi(r){\rm d}r<\infty \Longleftrightarrow \int_{[1,\infty)}r(\log r)\nu({\rm d}r)<\infty.$$
$$\int^1_0r^{-2}|\log r|\psi(r){\rm d}r<\infty\Longleftrightarrow\int_{[1,\infty)}r(\log r)^2\nu({\rm d}r)<\infty.$$
We can use these equivalences to provide some examples in which the moment conditions appearing in Theorem \ref{main-2} hold or fail.

Firstly, we provide an example where \eqref{A3}
 holds true but $\int_{[1,\infty)}r(\log r)\nu({\rm d}r)=\infty$. According to \cite{RW},
 $\psi_1(\lambda)=\lambda^2{\color{black}\log}^{-1}(1+\lambda),\lambda\ge 0$ is a branching mechanism.
 By some elementary calculations, we can check that $\psi_1$ satisfies \eqref{A3} but
 $\int^1_0r^{-2}\psi_1(r){\rm d}r=\infty.$ Let $\nu_1$ be the measure $\nu$ in \eqref{branch-mech}
 corresponding to $\psi_1$. Then
 $\int_{[1,\infty)}r\log r\nu_1({\rm d}r)=\infty$.

Secondly, we give an example where $\int_{[1,\infty)}r(\log r)\nu({\rm d}r)<\infty$ and $\int_{[1,\infty)}r(\log r)^2\nu({\rm d}r)=\infty$. According to \cite{RW},
$\psi_2(\lambda)=\lambda{\color{black}(\lambda\log\lambda-\lambda{\color{black}+1})/(\log \lambda)^2}, \lambda> 0$, is a branching mechanism.
We can check that $\psi_2$ satisfies \eqref{A3} and $\int^1_0r^{-2}\psi_2(r){\rm d}r<\infty$,
but $\int^1_0r^{-2}|\log r|\psi_2(r){\rm d}r=\infty$.  Let $\nu_2$ be the measure $\nu$ in \eqref{branch-mech}
corresponding to $\psi_2$.  Than $\int_{[1,\infty)}r\log r\nu_2({\rm d}r)<\infty$
but $\int_{[1,\infty)}r(\log r)^2\nu_2({\rm d}r)=\infty$.
\end{remark}

\section{Branching exit Markov systems and embedded continuous state branching processes}\label{BEM}

For each $y,t\geq 0$, define the space-time domain $D^t_y = \{(x,u) : x<y, \, u<t\}$
 and for each
$c\in\mathbb{R}$ let $X^c = \{X_t^c: t\geq 0\}$ be the
sequence of measures which satisfies $\l f, X_t^c\r =\l f(ct + \cdot), X_t \r$ for all $t\geq 0$ and $f\in C^+_b(\mathbb{R})$. It is straightforward to deduce that for each $\mu\in\mathcal{M}_F(\R)$,  $(X^c, \mathbb{P}_\mu)$ is a superprocess with general branching mechanism $\psi$ whose movement component corresponding to a Brownian motion with drift $c$.
According to Dynkin's theory of exit measures
\cite{Dyn2001} it is possible to describe the mass in the
superprocess $X^c$ as it first exits the growing family of domains
$\{D^t_y: t\geq 0, y\geq 0\}$ as a sequence of random measures on $\mathbb{R}\times[0,\infty)$, known as {\it branching Markov exit measures}, which we denote by $\{X^c_{D^t_y}
: t\geq 0, y\geq 0\}$. In particular, according to the characterisation  for branching Markov exit measures given in Section 1.1 of \cite{Dynkin-Kuznetsov}, each of the random measures $X^c_{D^t_y}$ is
 supported on $\partial D^t_y = (\{y\}\times[ 0,t)) \cup
((-\infty,y]\times\{t\})$ and has  the following defining Markov
branching property. Let $\mathcal{F}^c_{D^t_y} = \sigma(X^c_{D^u_x} :u\leq t, x\leq y)$.
For all $t\geq r$,  $y\geq z$, $\mu\in\mathcal{M}_F(\mathbb{R})$
{\color{black} with $\mbox{supp }\mu\subset (-\infty, z]$}
and
$f\in C_b^+(D^t_y)$,
\begin{equation}
\mathbb{E}_\mu (e^{-\langle f, X^c_{D^t_{y}}\rangle }|\mathcal{F}^c_{D^r_z})
= e^{-\langle u^y_f(\cdot, t-\cdot) , X^c_{D^r_z}\rangle},
\label{prePDE3}
\end{equation}
where, for all $(x,s)$ in  $ D^t_{y}$, $u^y_f$ is the unique positive solution of the partial differential equation
\begin{equation}
\frac{\partial}{\partial s}u^y_f(x,s)  = \dfrac{1}{2}\dfrac{\partial^2 }{\partial
       x^2}u^y_f(x,s) + c \frac{\partial}{\partial x}u^y_f(x,s)-\psi(u^y_f(x,s)),
       \label{PDE3a}
\end{equation}
with boundary conditions $u^y_f(y,s) =f(y, s)$ for $s\le t$   and $u^y_f(x,t) =f(x,t)$  for
$x\leq y$. We may similarly consider the branching Markov property
of the exit measures $X^c_{D^t_{-z}}$ where $D^t_{-z}=\{(x, r): r<
t, -z<x\}$ with $z\geq 0$. Moreover, define for convenience $D_y =
D^\infty_y$ and by monotonicity one may also define  $X^c_{D_y} =
\lim_{t\uparrow\infty}X^c_{D^t_y}|_{\{y\}\times [ 0,t)}$.

An important consequence of the Markov branching
property above is the following theorem
which will  feature crucially in our proof of Theorem
\ref{main-1}.

\begin{theorem}\label{CSBP} Define for each $y\geq 0$ and $c\geq 0$, $Z^c_y:=\langle 1, X^c_{D_y}\rangle = ||X^c_{D_y}||$ and $Z^c_{-y}:=\langle 1, X^c_{D_{-y}}\rangle = ||X^c_{D_{-y}}||$.  For all $x\in\mathbb{R}$ and $\lambda\in(0, \underline{\lambda}]$ the following statements hold $\mathbb{P}_{\delta_x}$-almost surely.
\begin{itemize}
\item[(i)] The process  $\{Z^{c_\lambda}_y: y\geq x\}$
 is a conservative supercritical continuous state
branching process with growth rate $\lambda$. Moreover, the process $Z^{c_\lambda}$ becomes extinct with positive probability if and only if (\ref{A3}) holds.
\item[(ii)] The process  $\{Z^{c_\lambda}_{y}: y\leq -x\}$  is a  subcritical continuous state branching process with growth rate $-\lambda$. Moreover, there is almost sure extinction if and only if  (\ref{A3}) holds.

\end{itemize}
\end{theorem}

\proof First part of (i). For $x\leq y$ and $f\in C^+_b(\R\times[0,\infty))$ such that $f(x,t) = f(x,0)=:f(x)$ for all $t\geq0$, let $v^y_f(x,t) := \mathbb{E}_{\delta_x}(\langle f, X^c_{D^t_y}\rangle).$ By performing a similar linearisation to the linearisation (\ref{PDE2}) of (\ref{PDE}), we have that
\begin{equation}
\frac{\partial}{\partial t} v^y_f(x,t) = \dfrac{1}{2}\dfrac{\partial^2 }{\partial
       x^2}v^y_f(x,t) + c \frac{\partial}{\partial x}v^y_f(x,t)-\psi'(0^+)v_f^y(x,t),
       \label{PDE3}
\end{equation}
with $v_f^y(y,s) = f(y)$ for
$s\le t$ and $v_f^y(x,0) = f(x)$ for $x\leq y$. The
classical Feynman-Kac formula allows us to write the unique solution
to (\ref{PDE3})  as
\begin{equation}\label{expect function}
v^{y}_f(x,t) = \Pi^c_x[e^{-\psi'(0^+)(t\wedge \tau^+_y)} f(\xi_{t\wedge \tau^+_y})],
\end{equation}
 where $\tau^+_y = \inf\{t> 0 : \xi_t >y\}$ and under $\Pi^c_x$,
$\{\xi_t : t\geq 0\}$ is a Brownian motion with drift $c$ issued
from $x$.  By means of an increasing sequence of continuous functions which are valued $0$ at $y$ and which converge pointwise to $\mathbf{1}_{(-\infty,y]}(\cdot)$, it is now possible to deduce by monotone convergence that
\begin{equation}
 \mathbb{E}_{\delta_x}(X^c_{D^t_y}((-\infty,y]\times\{t\})) = e^{-\psi'(0^+)t}\Pi_x^c(\tau^+_y>t).
\label{pre-y-in-t}
\end{equation}
For $x\leq y$ it is known that the density of
$\tau^+_y$ is given by
\begin{equation}
\frac{y-x}{\sqrt{2\pi t^3}}\exp\left(-\frac{(y-ct)^2}{2t}\right)
,\quad t>0.
\label{density}
\end{equation}
 Now let $c = c_\lambda$ for $\lambda \in(0,\underline\lambda]$. From (\ref{pre-y-in-t}) and (\ref{density}), an application of  {\color{black}L'H\^o}pital's rule shows that
 \begin{equation}
 \lim_{t\uparrow\infty}\mathbb{E}_{\delta_x}(X^{c_\lambda}_{D^t_y}((-\infty,y]\times\{t\}))  =
0.
\label{goes-to-zero-on-a-subsequence}
\end{equation}

It now follows from (\ref{expect function}) with $f=1$, (\ref{pre-y-in-t})  and (\ref{goes-to-zero-on-a-subsequence}) that
for all  $x\in(-\infty,y]$,
\[
\mathbb{E}_{\delta_x}(|| X^{c_\lambda}_{D_y}||) = \lim_{t\uparrow\infty}
\mathbb{E}_{\delta_x}(|| X^{c_\lambda}_{D^t_y}||) =\Pi^{c_\lambda}_x[e^{-\psi'(0^+)\tau^+_y}; \tau^+_y<\infty]= e^{\lambda (y-x)}.
\]
Note also from
(\ref{prePDE3}) we have that for all $a,b,\theta\geq0$ and
$x\in\mathbb{R}$,
\begin{equation}
\mathbb{E}_{(a+b)\delta_x}(e^{-\theta\langle 1, X^{c_\lambda}_{D_{y}}\rangle }) = e^{ -(a+b) v_\theta^y(x)} = \mathbb{E}_{a\delta_x}(e^{-\theta\langle 1, X^{c_\lambda}_{D_{y}}\rangle })\mathbb{E}_{b\delta_x}(e^{-\theta\langle 1, X^{c_\lambda}_{D_{y}}\rangle }),
\label{branching-property}
\end{equation}
 showing that $Z^{c_\lambda}$ is a conservative continuous state branching process.
\bigskip

First part of (ii). By symmetry, it suffices to prove that for
$\lambda\in(0, \underline{\lambda}]$, the process
$\{Z^{-c_\lambda}_{x}: x\geq 0\}$ is a subcritical
continuous state branching process with growth rate $-\lambda$. This
conclusion follows from a similar analysis to the proof of part (i),
noting in particular that
\[
\mathbb{E}_{\delta_x}(|| X^{-c_\lambda}_{D_y}||) =\Pi^{-c_\lambda}_x[e^{-\psi'(0^+)\tau^+_y}; \tau^+_y<\infty]= e^{-\lambda (y-x)}.
\]
The details are left to the reader.

\bigskip

Second part of (ii).
For any
$y\ge z$, $\mu\in {\cal M}(\mathbb{R})$ {\color{black} with $\mbox{supp } \mu\subset(-\infty, z]$}  and $\theta>0$,
\begin{equation}\label{pre-Lap-exit}
\mathbb{E}_\mu (e^{-\langle\theta, X^{-c_\lambda}_{D_{y}}\rangle }|\mathcal{F}^{-c_\lambda}_{D_z})  =  e^{-\langle u^y_\theta , X^{-c_\lambda}_{D_z}\rangle},
\end{equation}
where, by taking limits as $t$ and then $r$ tend to infinity in {\color{black}(\ref{PDE3a})}, we have that  $u^y_\theta$ is the unique positive solution to the equation
\begin{equation}\label{Lap-exit}
0= \dfrac{1}{2}\dfrac{\partial^2 }{\partial
       x^2}u^y_\theta(x) - c_\lambda \frac{\partial}{\partial x}u^y_\theta(x)-\psi(u^y_\theta(x)),
\end{equation}
on $(-\infty,y]$ with boundary value $u^y_\theta(y)=\theta$.

This tells us that for each fixed $\theta \geq 0$,
\[
\mathbb{E}(e^{-\langle\theta, X^{-c_\lambda}_{D_{x}}\rangle })  =
e^{-u^{x}_\theta (0)}  = e^{-u^0_\theta (-x)},
\]
where $u^0_\theta$ solves
\[
0= \dfrac{1}{2}\dfrac{\partial^2 }{\partial
       x^2}u^0_\theta(x)-c_\lambda \frac{\partial}{\partial x}u^0_\theta(x)-\psi(u^0_\theta(x)),
\]
on $(-\infty, 0)$ with boundary value
$u^0_\theta(0)=\theta$.
Written yet another way, this tells us that
$u^x_\theta : = u^x_\theta (0)$ satisfies
\[
0= \dfrac{1}{2}\dfrac{\partial^2 }{\partial
       x^2}u^x_\theta + c_\lambda \frac{\partial}{\partial x}u^x_\theta-\psi(u^x_\theta),
\]
on $(0,\infty)$ with boundary value $u^0_\theta = \theta$.  

On the other hand, if $\psi_{-c_\lambda}$ is the branching mechanism of $\{Z^{-c_\lambda}_{x} : x\geq 0\}$, then we also know that
\[
\frac{\partial}{\partial x}u^x_\theta + \psi_{-c_\lambda}(u^x_\theta)=0,
\]
for $x\geq 0$. Combining the previous two differential equations, we easily deduce that
\[
\frac{1}{2}\psi_{-c_\lambda}'(u^x_\theta)\psi_{-c_\lambda}(u^x_\theta)- c_\lambda \psi_{-c_\lambda}(u^x_\theta) = \psi(u^x_\theta).
\]
As $\{Z^{-c_\lambda}_{x} : x\geq 0\}$ is subcritical, we know that
$u^\infty_\theta = 0$. Thus by continuity, for each fixed
$\theta>0$, the range of $\{u^x_\theta :x \geq 0\}$ contains $[
0,\theta]$. As $\theta$ may be made arbitrarily large, it follows
that
\[
\frac{1}{4}\frac{\rm d}{{\rm d}u}\psi^2_{-c_\lambda}(u) - c_\lambda \psi_{-c_\lambda}(u) = \psi(u), \, \, u\geq 0.
\]
Subcriticality also implies that  $\psi_{-c_\lambda}(0) = 0$.

Next note that
\begin{equation}
\psi^2_{-c_\lambda}(u)- c_\lambda \int_{\lambda^*}^u \psi_{-c_\lambda}(s){\rm d}s = \int_{\lambda^*}^u \psi(s){\rm d}s.
\label{dividebyintegral}
\end{equation}
As $\psi_{-c_\lambda}$  tends to infinity at infinity, we may apply L'H\^opital's rule to deduce that
\begin{equation}
\lim_{u\uparrow\infty}\frac{\int_{\lambda^*}^u \psi_{-c_\lambda}(s){\rm d}s}{{\color{black}\psi^2_{-c_\lambda}(u)}}  = \lim_{u\uparrow\infty}\frac{1}{\psi'_{-c_\lambda}(u)}.
\label{ratiolimit}
\end{equation}
Note it follows in a straightforward manner form the L\'evy-Khintchine formula
that the limit on the right hand side above exists (and may possibly equal zero).
The limit (\ref{ratiolimit}) when combined with (\ref{dividebyintegral}) now allows us to conclude that
\begin{equation}
\int^\infty \frac{1}{\psi_{-c_\lambda}(\xi)}{\rm d}\xi<\infty
\Longleftrightarrow
\int^\infty \frac{1}{\sqrt{\int_{\lambda^*}^\xi\psi(u){\rm d}u}}{\rm d}\xi<\infty.
\label{replace}
\end{equation}
As $\{Z^{-c_\lambda}_{x} : x\geq 0\}$ is subcritical, this is equivalent to saying that there is almost sure extinction if and only if (\ref{A3}) holds.

\bigskip

Second part of (i). Using exactly the same proof we can show that (\ref{replace}) holds with $\psi_{-c_\lambda}$ replaced by $\psi_{c_\lambda}$. The desired result follows by recalling that $\int^\infty 1/ \psi_{c_\lambda}(\xi){\rm d}\xi<\infty$ is the necessary and sufficient condition in the current context for the event of extinction to agree with the event of becoming extinguished.
\qed

\begin{cor}\label{support} Suppose that (\ref{A3}) holds. Fix $x\in\mathbb{R}$.
For each $c\geq 0$, let $\mathcal{R}^{-c}$ be the smallest closed set containing ${\rm supp}X^{-c}_t$ for all $t\geq 0$.
Then for all $c\geq
c_{\underline\lambda}=\underline\lambda$ we have that
$\mathbb{P}_{\delta_x}(\sup\mathcal{R}^{-c}<\infty) = 1$ and for all $c
<c_{\underline\lambda}$ we have
$\mathbb{P}_{\delta_x}(\sup\mathcal{R}^{-c}=\infty| \mathcal{E}^c) = 1$.
In particular if $R_t = \sup\{y\in\mathbb{R}: X_t(y, \infty)>0\}$
then
\begin{equation}
\lim_{t\uparrow\infty}\frac{R_t}{t} = \underline\lambda
\label{sym1},
\end{equation}
$\mathbb{P}_{\delta_x}$-almost surely  on  $\mathcal{E}^c$.
\end{cor}

\proof  From  Theorem \ref{CSBP} (ii), under the assumption of (\ref{A3}), the
process $\{Z^{-c_\lambda}_{x}: x\geq 0\}$ is subcritical and
becomes extinct with probability 1. This implies that for all $c\geq
c_{\underline\lambda}=\underline\lambda$,
$\mathbb{P}_{\delta_x}(\sup\mathcal{R}^{-c}<\infty) = 1$ and hence
\[
\limsup_{t\uparrow\infty}\frac{R_t}{t}\leq \underline\lambda,
\]
where $R_t = \sup\{y\in\mathbb{R}: X_t
(y,\infty) >0\}$.

Next we want to show
\begin{equation}
\liminf_{t\uparrow\infty}\frac{R_t}{t} \geq \underline\lambda,
\label{liminfR}
\end{equation}
on $\mathcal{E}^c$.

We shall use the conclusion of Theorem \ref{mgcgce} to prove (\ref{liminfR}). The reader should note that the proof of Theorem \ref{mgcgce}, which appears later in this paper, does not depend on the result we are currently proving.
We also make use of an argument which is essentially taken from Git et al. \cite{GHH}.
For $0<\epsilon < \underline\lambda /2$ and  $\gamma = \underline\lambda -
\epsilon$ note that $
e^{\gamma x }\mathbf{1}_{(x\leq (\gamma - \epsilon)t)} \leq
e^{(\gamma-\epsilon) x }
e^{\epsilon(\gamma-\epsilon)t},
$ and hence
\begin{equation}\label{extreme}
\limsup_{t\uparrow\infty}
e^{-(\gamma^2 /2 -\psi'(0^+))t}
\langle
e^{\gamma\cdot}\mathbf{1}_{(\cdot\leq (\gamma -\epsilon)t)} , X_t
\rangle
\leq \limsup_{t\uparrow\infty}
e^{-\epsilon^2 t/2}
W_t (-\gamma+\epsilon) = 0,
\end{equation}
$\mathbb{P}_{\delta_x}$-almost surely.
It follows that
\begin{equation}\label{extreme-II}
\lim_{t\uparrow\infty}
e^{-(\gamma^2 /2 - \psi'(0^+))t}
\langle
e^{\gamma\cdot}\mathbf{1}_{(\cdot > (\gamma -\epsilon)t )} , X_t
\rangle
= W_\infty(-\gamma+\epsilon),
\end{equation}
$\mathbb{P}_{\delta_x}$-almost surely. Note that by Theorem \ref{mgcgce} (i)  the event $\{W_\infty(-\gamma+\epsilon)>0\}$ agrees with $\mathcal{E}^c$. Hence as $\epsilon$ can be made arbitrarily small, (\ref{liminfR}) follows on $\mathcal{E}^c$ .

Together with (\ref{liminfR}) this implies the strong law of large numbers, (\ref{sym1}), on $\mathcal{E}^c$ and all other claims in the corollary follow immediately. \qed

\begin{theorem}\label{same-survival}
Suppose that $\lambda\in(0,\underline{\lambda}]$ {\color{black}and $Z^{c_\lambda}=\{Z^{c_\lambda}_y:y\ge 0\}$.}
Then $\mathbb{P}$-almost surely,  $\{Z^{c_\lambda} \text{ extinguishes}\}$
agrees with the event $\mathcal{E}$.
\end{theorem}

\proof
First we establish that $\widetilde{\mathcal{E}}:=\{Z^{c_\lambda}\text{ extinguishes}\}\subseteq\mathcal{E}$. Begin by noting that, thanks to  monotonicity,
\begin{equation}
\lim_{y\uparrow\infty}X^{c_\lambda}_{D^t_y}((-\infty, y]\times\{t\}) = ||X_t||,
\label{thirdline}
\end{equation}
Next note that since $Z^{c_\lambda}$ is a supercritical conservative branching process, it follows that there is a $\lambda_0>0$ such that $\mathbb{P}(\widetilde{\mathcal{E}}) = e^{-\lambda_0}$. Using the Markov branching property for exit measures we have,
\[
\mathbb{E}(\mathbf{1}_{\widetilde{\mathcal{E}}} | \mathcal{F}^{c_\lambda}_{D^t_y}) =
e^{-\lambda_0 || X_{D^t_y}^{c_\lambda}||}\leq e^{-\lambda_0 X_{D^t_y}^{c_\lambda}((-\infty,y]\times\{t\}) }.
\]
Hence
\begin{eqnarray}
\mathbb{E}\left[\lim_{t\uparrow\infty}\lim_{y\uparrow\infty}\mathbb{E}(\mathbf{1}_{\widetilde{\mathcal{E}}} | \mathcal{F}^{c_\lambda}_{D^t_y})\mathbf{1}_{\mathcal{E}^c} \right]&\leq& \mathbb{E}\left[\lim_{t\uparrow\infty}\lim_{y\uparrow\infty} e^{-\lambda_0 X_{D^t_y}^{c_\lambda}((-\infty,y]\times\{t\}) }\mathbf{1}_{\mathcal{E}^c}\right]\notag\\
& =& \mathbb{E}\left[\lim_{t\uparrow\infty}e^{-\lambda_0 ||X_t|| }\mathbf{1}_{\mathcal{E}^c}\right]\notag\\
&=&0.\label{inlightof}
\end{eqnarray}
Note that in the first equality we have used the fact that $||X^{c_\lambda}_t|| = ||X_t||$.
Our objective is to show that
\[
\mathbb{E}(\mathbf{1}_{\mathcal{E}^c \cap\widetilde{\mathcal{E}}}) = 0,
\]
and hence $\widetilde{\mathcal{E}}\subseteq\mathcal{E}$, $\mathbb{P}$-almost surely. To this end, in light of (\ref{inlightof}), it suffices to prove that $\mathcal{E}^c\in\sigma\left(\bigcup_{t>0}\bigcup_{y>0}\mathcal{F}^{c_\lambda}_{D_y^t}\right)$.
Note however that,  by (\ref{thirdline}),   $||X_t||\in\sigma\left(\bigcup_{y>0}\mathcal{F}^{c_\lambda}_{D_y^t}\right)$, which implies that
 ${\cal E}\in\sigma\left(\bigcup_{t>0}\bigcup_{y>0}\mathcal{F}^{c_\lambda}_{D_y^t} \right)$.

\medskip

Now fix $t>0$. Note that the Markov branching property applied to the exit measure  $X^{c_\lambda}_{D^t_y}$ implies that
\[
 \mathbb{P}(\mathcal{E}) = \mathbb{E}(\mathbb{P}(\mathcal{E}|\mathcal{F}^{c_\lambda}_{D^t_y})) = \mathbb{E}\left[e^{-\lambda^*||X^{c_\lambda}_{D^t_y}||}\right]
 \leq \mathbb{E}\left[e^{-\lambda^*X^{c_\lambda}_{D^t_y}(\{y\}\times[ 0,t))}\right].
 \]
Now taking limits as $t\uparrow\infty$ we have with the help of both monotone and dominated convergence that
\[
\mathbb{P}(\mathcal{E})\leq \mathbb{E}\left[e^{-\lambda^*||X^{c_\lambda}_{D_y}||}\right] = \mathbb{E}(e^{-\lambda^*Z^{c_\lambda}_y} ).
\]
Taking limits again as $y\uparrow\infty$ we find that $\mathbb{P}(\mathcal{E})\leq \mathbb{P}(Z^{c_\lambda}\text{ extinguishes})$. Together with the conclusion of the previous paragraph, we are forced to conclude that $\mathcal{E} = \{Z^{c_\lambda} \text{ extinguishes}\}$, $\mathbb{P}$-almost surely, as required.
\qed

\section{Proof of Theorem \ref{main-1}}

\noindent{\bf Proof of Theorem \ref{main-1} (i):} Suppose that there
exists a travelling wave at speed $c\in[ 0,c_{\underline\lambda})$
which we shall denote by $\Phi$. 
For all $x\in\mathbb{R}$ and $t\geq 0$, we have
$
 \mathbb{E}_{\delta_x}(e^{-\langle\Phi, X^c_t\rangle}) = e^{-u^c_\Phi(x,t)},
$
where $u^c_\Phi$ solves
\begin{equation}
\frac{\partial}{\partial t}u^c_\Phi(x,t) = \frac{1}{2}\frac{\partial^2}{\partial x^2}u_\Phi^c(x, t)+ c\frac{\partial}{\partial x}u_\Phi^c(x,t) -\psi(u^c_\Phi(x,t)),
\label{PDEdrift}
\end{equation}
with initial condition $u(x,0) = \Phi(x)$.
This partial differential equation has a unique positive solution for the same  reasons that (\ref{PDE}) has a unique solution. Since $\Phi(x)$ also solves (\ref{PDEdrift}), it follows that
\[
 \mathbb{E}_{\delta_x}(e^{-\langle\Phi, X^c_t\rangle}) = e^{-\Phi(x)}.
\]
Together with the branching property, this in turn implies that $\{e^{-\langle\Phi, X^c_t\rangle}:t\geq 0\}$ is a uniformly integrable martingale. Its   almost sure and $L^1(\mathbb{P}_{\delta_x})$ limit is denoted by $M_\infty$.
The Markov branching property applied to the exit measure $X^c_{D^t_y}$ implies that for all $x\leq y$,
\[
\mathbb{E}_{\delta_x}\left(e^{-\langle\Phi, X^c_{D^t_y}\rangle} \right)= \mathbb{E}_{\delta_x}\left[\mathbb{E}(M_\infty| \mathcal{F}^c_{D_y^t})\right] = e^{-\Phi(x)}.
\]
Note however that for all $z\in{\rm supp} X^c_{D^t_y}$ we have by
monotonicity, $\Phi(z)\geq \Phi(y)$. Moreover, as a measure, we also have  $X^c_{D^t_y}\geq
X^c_{D^t_y}|_{{\color{black}(-\infty,y)}\times\{t\}}$.  It follows that
\begin{equation}
e^{-\Phi(x)}\leq
\mathbb{E}_{\delta_x}\left(e^{-\Phi(y)X^c_{D^t_y}({\color{black}(-\infty,y)}\times\{t\})}
\right). \label{ytoinf}
\end{equation}
{\color{black}
Our next objective is to show that for any $y>x$,
\begin{equation}\label{toshow}\mathbb{P}_{\delta_x}(\liminf_{n\uparrow\infty}X^c_{D^{n}_y}((-\infty,y)\times\{n\})  >0)>0 .
\end{equation}
Suppose now that the probability in (\ref{toshow}) is equal to zero for a given $y>x$. Let $X(n): =X^c_{D^{n}_y}((-\infty,y)\times\{n\})$ for $n\geq 0$. Then
\begin{equation}\label{liminf-0}\liminf_{n\to\infty}X(n)=0, \quad \mathbb{P}_{\delta_x}\mbox{-a.s.}\end{equation} Note that, under (5),  $0$ is an absorbing state for the sequence $\{X(n): n\geq 0\}$ in the sense that $X(m) = 0$ implies that $X(m+k) = 0$ for all $k\geq 0$. Note  that since $\{X^c_{D^t_y}|_{{\color{black}(-\infty,y)}\times\{t\}}: t\geq 0\}$ is a superprocess with branching mechanism $\psi$ and underlying motion which is that of a Brownian motion with drift $c$ killed on hitting $y$, and therefore Markovian, then we have the estimate
\begin{eqnarray*}
\lefteqn{\mathbb{P}_{\delta_x}(\exists~ m \text{ s.t. } X(m) = 0| X(0),\ldots, X(n)) }&&\\
&&\geq \inf_{\mu: ||\mu|| = X(n)}\mathbb{P}_\mu(\exists~ m\text{ s.t. } X(m) = 0)\\
&&\geq \inf_{\mu: ||\mu|| = X(n)}\mathbb{P}_\mu(\mathcal{E})\\
&& = e^{-\lambda^* X(n)}.
\end{eqnarray*}
  Letting $n\to\infty$ in the above inequality, by \eqref{liminf-0}, we obtain that $P_{\delta_x}(\exists~ m \text{ s.t. } X(m) = 0)=1$. It follows that $\liminf_{t\uparrow\infty}L^c_t \ge y$ $\mathbb{P}_{\delta_x}$-almost surely where
$L^c_t = \inf\{z: X^c_t(-\infty,z] >0\}$. However, from Corollary \ref{support}, we also deduce that under (\ref{A3}),
\[
\lim_{t\uparrow\infty}\frac{L^c_t}{t} = c- \underline{\lambda}<0,
\]
which constitutes a contradiction. Therefore \eqref{toshow} holds for any $y>x$.

It follows by \eqref{toshow} and
the Reverse Fatou Lemma, that
}
\[
e^{-\Phi(x)}\leq \limsup_{n\uparrow\infty}
\mathbb{E}_{\delta_x}\left(e^{-\Phi(y)X^c_{D^{n}_y}((-\infty,y]\times\{n\})}
\right) <1,
\]
for all sufficiently large $y>x$. 
As $x$ may be chosen arbitrarily in this argument, it follows that there exists a constant $C>0$ such that $\Phi(x)>C$ for all $x\in\mathbb{R}$. This leads to a contradiction of the assumption that $\Phi$ is a travelling wave.
\qed
\bigskip

\noindent{\bf Proof of Theorem \ref{main-1} (ii):} Grey \cite{Grey}
solves the classical Seneta-Heyde norming problem for continuous
state branching processes. In particular he shows that for all
$\lambda\in(0, \underline{\lambda}]$, taking account
of the fact that $Z^{c_\lambda}$ is a continuous state branching
process with growth rate $\lambda$, there exists a slowly varying
function at $0$, $L_\lambda$ such that
\begin{equation}
\lim_{x\uparrow\infty}e^{-\lambda x}L_\lambda(e^{-\lambda x})Z^{c_\lambda}_x = \Delta(\lambda),
\label{SHnorming}
\end{equation}
where $\Delta(\lambda)\geq 0$ is non-degenerate and the event
$\{\Delta(\lambda) = 0\}$ agrees with the event that $Z^{c_\lambda}$
becomes extinguished which in turn, by Theorem \ref{same-survival} agrees with the event $\mathcal{E}$.  
Note from  (\ref{SHnorming}) and the
fact that $L_\lambda$ is slowly varying, it is straightforward to
show that for all $\mu\in\mathcal{M}_F(\mathbb{R})$
\[
\mathbb{E}_\mu\bigg[\exp\{-\Delta(\lambda)\}\bigg] = \exp\{ - \langle \Phi,\mu\rangle\},
\]
where for all $x\in\mathbb{R}$,
\begin{equation}\label{auxiliar}
e^{-\Phi(x)}=
\mathbb{E}_{\delta_x}[\exp\left\{-\Delta(\lambda)\right\}]
 = \mathbb{E}[\exp\left\{-e^{-\lambda x}\Delta(\lambda)\right\}].
\end{equation}

Note in particular that $\Phi$ is a monotone decreasing function which is twice continuously differentiable on $(0,\infty)$ and moreover satisfies $\Phi(\infty) = 0$ and $\Phi(-\infty) = \lambda^*$.

The Markov branching property for $Z^{c_\lambda}$ implies that  by conditioning on $\mathcal{F}^{c_{\lambda}}_{D_x}$, where $x,z\in\mathbb{R}$, we get
\begin{equation}\label{martingale-first passage}
e^{-\Phi(z)}=\mathbb{E}\bigg[\exp\left\{-e^{-\lambda
z}\Delta(\lambda)\right\}
\bigg]
=\mathbb{E}\bigg[
\exp\{-\Phi(x+z)Z^{c_\lambda}_x\}
\bigg].
 \end{equation}

Setting $z=0$, $\mu=\delta_0$, $f = \Phi$, $c=c_\lambda$ {\color{black} in \eqref{pre-Lap-exit} and \eqref{Lap-exit} we see  that $\Phi$ necessarily solves}
\begin{equation}\label{kpp-equation}
\frac{1}{2}\Phi''+c_\lambda\Phi'-\psi(\Phi)=0,\, \mbox{ on }\mathbb{R},
\end{equation}
as required.
\qed

\noindent{\bf Proof of Theorem \ref{main-1} (iii) and (iv):} Let
$\lambda\in (0, \underline\lambda]$ and assume that
$\Phi_{c_\lambda}$ is a
  travelling wave solution to the equation \eqref{travwave} with speed
$c_\lambda.$  From (\ref{prePDE}) and (\ref{PDE}) it follows that  for all $z\in\mathbb{R}$.
\[
e^{-\Phi_{c_\lambda}(z)}  = \mathbb{E}_{\delta_z} [\exp\{ -\langle \Phi_{c_\lambda}(\cdot + c_\lambda t), X_t\rangle\}] =
\mathbb{E}_{\delta_z} [\exp\{ -\langle \Phi_{c_\lambda}, X^{c_\lambda}_t\rangle\}]
= \mathbb{E}[\exp\{ -\langle \Phi_{c_\lambda}(z+\cdot), X^{c_\lambda}_t\rangle\}].
\]
and hence, together with the branching Markov property, we have that for all $z\in\mathbb{R}$
\begin{equation}
M_t^{\lambda, z}:  =\exp\{ - \langle \Phi_{c_\lambda}(z+\cdot), X^{c_\lambda}_t\rangle\}, \, t\geq 0
\label{anothermg}
\end{equation}
is a positive uniformly integrable $\mathbb{P}$-martingale.
From (\ref{goes-to-zero-on-a-subsequence}) we may deduce that there exists a deterministic  subsequence $\{t_n : n\geq 0\}$ which increases to infinity (and may depend on $x$) such that
\begin{equation}\label{extinguish}
\lim_{n\uparrow\infty}X^{c_\lambda}_{D^{t_n}_x}((-\infty, x]\times\{t_n\})=0.
\end{equation}
Now fix $z\in\mathbb{R}$ and $x\geq 0$. Let $M_\infty^{\lambda, z} :=
\lim_{t\uparrow\infty} M_t^{\lambda, z}$ then the branching Markov
property applied to the exit measure
$X_{D_z^{t_n}}^{c_\lambda}$ gives us
\[
\mathbb{E}\left[M_\infty^{\lambda,
z}\bigg|\mathcal{F}^{c_\lambda}_{D_x}\right]
=\lim_{n\uparrow\infty}\mathbb{E}\left[M_\infty^{\lambda,
z}\bigg|\mathcal{F}^{c_\lambda}_{D^{t_n}_x}\right]
=\lim_{n\uparrow\infty}\exp\{-\langle\Phi_{c_\lambda},
X^{c_\lambda}_{D^{t_n}_x}\rangle\}
 = \exp\{-\Phi_{c_\lambda}(x+z)Z^{c_\lambda}_x\}.
 \]
 From (\ref{anothermg}) we see that the event $\{M^{\lambda,z}_\infty=1\}$ contains in the event
 that $X$ becomes extinguished which in turn, from the proof of the previous part of the theorem,  agrees with the event that
 $\{\Delta(\lambda)=0\}$. Recalling (\ref{SHnorming}),
  it follows that there exist a set of positive $\mathbb{P}$-probability on which $Z^{c_\lambda}$ has a strictly positive normalised limit, such that the normalising sequence may be taken to be either $e^{-\lambda x}L_\lambda(e^{-\lambda x})$ or $\Phi_\lambda(x+z)$ as $x\uparrow\infty$.
 It must therefore follow that there exists a constant $k_{\lambda,z}\in(0,\infty)$ such that
 \[
\lim_{x\uparrow\infty}\frac{ \Phi_\lambda(x+z)}{e^{-\lambda x} L_\lambda (e^{-\lambda x})} = k_{\lambda,z}.
 \]
 As $L_\lambda$ is slowly varying it is easy to see that $k_{\lambda, z} = e^{-\lambda z}k_\lambda$ where $k_\lambda := k_{\lambda, 0}$.

 This also tells us that
 \[
 M^{\lambda,z}_\infty = \exp\{ - k_\lambda e^{-\lambda z}\Delta(\lambda)\},
 \]
 and hence taking expectations with respect to $\mathbb{P}$ we see that
 \[
e^{- \Phi_{c_\lambda}(z)} = \mathbb{E}\bigg[\exp\{ - k_\lambda e^{-\lambda z}\Delta(\lambda) \}\bigg],
 \]
 thus establishing uniqueness up to an additive constant.
\qed

Reviewing  the proof of Theorem \ref{main-1} we obtain the following corollary, a simpler version of which has appeared in Neveu \cite{neveu} and in parallel to writing of this paper, a similar result for branching Brownian motion has been described in \cite{Mal}.
\begin{cor}
For $\lambda\in(0,\underline{\lambda}]$, the
continuous state branching process $Z^{c_\lambda}$ has branching
mechanism
\[
\psi_{c_\lambda}(\theta) = \Phi_{c_\lambda}'(\Phi^{-1}_{c_\lambda}(\theta)),
\]
for $\theta\in[ 0,\lambda^*]$ where $\Phi_{c_\lambda}$ is any
version of the unique travelling wave at speed $c_\lambda$.
Alternatively $\psi_{c_\lambda}$solves the differential equation
\[
\dfrac{1}{4}\frac{\rm d}{{\rm d}u}f^2(u)+c_\lambda f(u)=\psi(u),\qquad
u\in(0, \lambda^*).
\]
with boundary conditions $f(0) = 0$ and $f(\lambda^*) = 0$.
\end{cor}
\proof
Since  $Z^{c_\lambda}$ is a continuous time
continuous state branching process with branching mechanism, say $\psi_{c_\lambda}$, equation (\ref{martingale-first passage}) implies that
$\Phi(z)=u_x(\Phi(x+z)),$ where for $\theta\geq 0$,  $u_x(\theta)$ satisfies the semi-group equation
\begin{eqnarray}\label{cont dual equation}
\dfrac{\partial u_x(\theta)}{\partial x}+\psi_{c_\lambda}(u_x(\theta))=0,
\end{eqnarray}
with initial condition $u_0(\theta) = \theta$.
Differentiating the equality $\Phi(z) = u_x(\Phi(x+z))$ with respect to $x$ we get
$$0 =
\dfrac{\partial u_x(\Phi(z+x))}{\partial x}+\dfrac{\partial
u_x(\theta)}{\partial \theta}\bigg|_{\theta=\Phi(z+x)}\Phi'(z+x).
$$
By setting $x=0$ in the previous equality, making use of (\ref{cont dual equation}) and the fact that
 $\partial u_0(\theta)/\partial \theta=1$,
 we obtain that
\begin{equation}\label{first order der}
\Phi'(z)=\psi_{c_\lambda}(\Phi(z)).
\end{equation}

The first part follows directly from (\ref{first order der}). For the second part, one may differentiate (\ref{first order der}) and obtain
\begin{equation}\label{second order der}
\Phi''(z)=\psi'_{c_\lambda}(\Phi(z))\psi_{c_\lambda}(\Phi(z)).
\end{equation}
Combining \eqref{first order der} and \eqref{second order der} with
\eqref{kpp-equation} and noting the domain of $u:=\Phi(x), x\in \mathbb{R}$
is $(0, \lambda^*)$, we get
\[
\dfrac{1}{2}\psi_{c_\lambda}'(u)\psi_{c_\lambda}(u)+c_\lambda\psi_{c_\lambda}(u)-\psi(u)=0,
\]
which is the claimed differential equation.

To show the boundary conditions, first note that in all cases,
since, by Theorem \ref{CSBP},  $Z^{c_\lambda}$ is conservative, we
necessarily have $\psi_{c_\lambda}(0) = 0$. For the other boundary
condition recall from the beginning and the end of the proof of
Theorem \ref{main-1} (ii) that the event $\{\Delta(\lambda) = 0\}$
agrees both with the event $\mathcal{E}$  as well as the event that
$Z^{c_\lambda}$ becomes extinguished. Hence we have $\mathbb{P}(
Z^{c_\lambda} \text{ becomes extinguished}) =e^{-\lambda^*}$ and
then necessarily $\lambda^*$ must be the largest root in $[
0,\infty)$ of the equation $\psi_{c_\lambda}(\theta) = 0$. This
justifies the boundary condition $\psi_{c_\lambda} (\lambda^*) = 0$.
\qed

\section{Pathwise  spine decomposition}\label{spinesection}

The convergence of the martingales $W(\lambda)$ and $\partial
W(\lambda)$ to non-trivial limits will ultimately allow us to
identify the limiting variable $\Delta(\lambda)$ as either
$W_\infty(\lambda)$ or $\partial W_\infty(\lambda)$. There is a well
understood technique for branching particle processes, due to Lyons
et al. \cite{LPP} and Lyons \cite{Ly}, which can be used to
establish in a straightforward way conditions under which the latter
limits are non-trivial. This involves looking at how the given
martingales (or variants of them) perform as changes of measure. In
the case of superprocesses, this technique can be seen under the
pretext of Evans' immortal particle decomposition; see for example
Evans \cite{Evans} and Engl\"ander and Kyprianou \cite{EK}. In the
latter references, the decomposition has only been explored through
the semi-group of the underlying superprocess which has its
limitations when using it to analyse the martingales $W(\lambda)$
and $\partial W(\lambda)$ in the spirit of pathwise spine
decompositions for branching particle processes. In the analysis
below, we also use a new immortal particle decomposition for our
class of superprocesses which is defined in a pathwise sense and
therefore lends itself to the aforementioned classical martingale
analysis.  The feature which is in particular new to our spine
decomposition is the use of the Dynkin-Kuznetsov
$\mathbb{N}$-measure to describe a Poisson point process of
immigration along the immortal particle. For  branching mechanism without quadratic term (i.e,  $\beta=0$ in the definition of $\psi$ given by \eqref{branch-mech}), a similar pathwise  spine decomposition was given by Liu et al. \cite{Rong-Liet.al.} when dealing with another martingale which, like $W(\lambda)$, was constructed from a positive eigen-function to the linear semi-group of the underlying superprocess.

\bigskip

Let us move to our new spine decomposition and hencewith we start by defining some martingale
changes of measure. For each $\lambda\in\mathbb{R}$ and $\mu\in \mathcal{M}_F(\mathbb{R})$,
let $\mathbb{P}^{-\lambda}_{\mu}$ be defined by
\begin{equation*}
 \frac{{\rm d}\mathbb{P}^{-\lambda}_{\mu}}{{\rm d}\mathbb{P}_{\mu}}\bigg{|}_{\F_t}=\frac{W_t(\lambda)}{W_0(\lambda)},\,\, t\geq 0,
\end{equation*}
where $\mathcal{F}_t = \sigma(X_s : s\leq t)$. Note in particular that $W_0(\lambda) = \langle e^{-\lambda\cdot}, \mu\rangle$.
Next recall that for each $x\in\mathbb{R}$ we defined the process $\xi:=\{\xi_t : t\geq 0\}$ under $\Pi_x$ to be a Brownian motion issued from $x$.
If $\Pi_x^{-\lambda}$ is the law under which $\xi$ is a Brownian  motion with drift $-\lambda\in\mathbb{R}$ issued from $x\in\mathbb{R}$, then
 for each $t\geq 0$,
\begin{equation}
  \frac{{\rm d}\Pi^{-\lambda}_x}{{\rm d}\Pi_x}\bigg{|}_{\mathcal{G}_t}=e^{-\lambda
  (\xi_t-x)-\lambda^2 t/2},
  \label{girsanov}
\end{equation}
where  $\mathcal{G}_t=\sigma(\xi_s,s\leq t)$.
For convenience we shall also write
\begin{equation}
\Pi^{-\lambda}_\mu (\cdot) = \frac{1}{\langle e^{-\lambda \cdot}, \mu\rangle}\int_{\mathbb{R}} e^{-\lambda x}\mu({\rm d}x)\Pi^{-\lambda}_x(\cdot).
\label{randomization}
\end{equation}

              \begin{teo} \label{spineth}
              Suppose that $\lambda \in\mathbb{R}$,
              $\mu\in\M_F(\R)$ and
              $g\in C_b^+(\R)$.  Then
                          \begin{equation}\label{spine}
                           \mathbb{E}^{-\lambda}_{\mu}\left[e^{-\l g,X_t\r}\right]
                          =\mathbb{E}_{\mu}\left[e^{-\l g,X_t\r}\right]
                         \Pi_\mu^{-\lambda}\left[\exp\left\{-\int_0^t\phi(u_g(\xi_{t-s},s)){\rm d}s\right\}\right],
                          \end{equation}
                          where
\begin{equation}\label{phi}
\phi(\lambda) = \psi'(\lambda) - \psi'(0^+) =
2\beta\lambda + \int_0^\infty (1- e^{-\lambda x})x\nu({\rm d}x).
 \end{equation}
 and $u_g$ is the unique solution of (\ref{PDE2}).
               \end{teo}
\proof
       By the definition of $\mathbb{E}^{-\lambda}_{\mu}$, we get
       \begin{eqnarray*}
        \mathbb{E}^{-\lambda}_{\mu}\left(e^{-\l g,X_t\r}\right)
    &=& \frac{1}{\langle e^{-\lambda\cdot}, \mu\rangle}\mathbb{E}_{\mu}\left(e^{-\lambda c_\lambda t}\l e^{-\lambda\cdot},X_t\r e^{-\l g,X_t\r}\right)\\
    &=&-\frac{1}{\langle e^{-\lambda\cdot}, \mu\rangle} e^{-\lambda c_\lambda t} \mathbb{E}_{\mu}\left(\frac{\partial}{\partial\theta}e^{-\l g_\theta,X_t\r}
          \bigg|_{\theta=0^+}\right),
       \end{eqnarray*}
with  $g_\theta(x):=g(x)+\theta e^{-\lambda x}$.
Interchanging the expectation and differentiation, we get
       \begin{eqnarray*}
              \mathbb{E}^{-\lambda}_{\mu}\left(e^{-\l g,X_t\r}\right)
           &=&-\frac{1}{\langle e^{-\lambda \cdot}, \mu\rangle}{e^{-\lambda c_\lambda t}}\frac{\partial}{\partial\theta}
           e^{- \langle u_{g_\theta}(\cdot,t) , \mu\rangle}\big|_{\theta=0^+},
\end{eqnarray*}
where $u_{g_\theta}$ satisfies (\ref{PDE}) with $g$ replaced by
$g_\theta$.
Note that, $u_{g_0}=u_{g}$.  Hence,
           \begin{eqnarray}
                 \mathbb{E}^{-\lambda}_{\mu}\left(e^{-\l g,X_t\r}\right)
              &=&\frac{1}{\langle e^{-\lambda \cdot}, \mu\rangle}{e^{-\lambda c_\lambda t}}e^{- \langle u_{g_\theta}(\cdot,t), \mu\rangle}
                 \frac{\partial}{\partial\theta}\langle u_{g_\theta}(\cdot,t), \mu\rangle\big|_{\theta=0^+}
     \label{put-m-in}.
           \end{eqnarray}

Now let $m^g(x,t):=\dfrac{\partial}{\partial\theta}u_{g_\theta}(x,t)\big|_{\theta=0^+}$ for all $x\in\mathbb{R}$.
Taking derivatives in (\ref{PDE}) with $g$ replaced by $g_\theta$ and then taking the limit as $\theta$ goes to zero, we obtain the differential equation
             \begin{eqnarray*}
             \begin{cases}
             \dfrac{\partial}{\partial t}m^g(x,t)
         =\dfrac{1}{2}\dfrac{\partial^2}{\partial
         x^2}m^g(x,t)-\psi^\prime(u_g(x,t))m^g(x,t),\\
         m^g(x,0)=e^{-\lambda x}.
             \end{cases}
             \end{eqnarray*}
The classical Feynman-Kac formula gives
             \begin{equation*}
            m^g(x,t)=
            \Pi_x\left[e^{-\lambda \xi_t}\exp\left\{-\int_0^t\psi^\prime(u_g(\xi_{t-s},s)){\rm d}s\right\}\right].
             \end{equation*}
Plugging back into (\ref{put-m-in})  yields the following equality,
        \begin{eqnarray*}
      \lefteqn{   \mathbb{E}^{-\lambda}_{\mu}\left(e^{-\l g,X_t\r}\right) }&&\\
        &=&\frac{1}{\langle e^{-\lambda\cdot}, \mu\rangle}\mathbb{E}_{\mu}\left[e^{-\l g,X_t\r}\right]
\int_{\mathbb{R}} e^{-\lambda x} \mu({\rm d}x)        \Pi_x\left[{e^{-\lambda( \xi_t-x)-\lambda c_\lambda t}}
          \exp\left\{-\int_0^t\psi^\prime(u_g(\xi_{t-s},s)){\rm d}s\right\}\right]\\
       &=&   \mathbb{E}_{\delta_x}\left[e^{-\l g,X_t\r}\right]
         \Pi_\mu^{-\lambda}\left[
          \exp\left\{-\int_0^t\phi(u_g(\xi_{t-s},s)){\rm d}s\right\}\right],
        \end{eqnarray*}
        where in the final equality we have used (\ref{girsanov}) and the fact that  $\lambda
c_\lambda=-\psi^\prime(0^+)+\lambda^2/2$.

To complete the proof, note that the expression given for $\phi(\lambda)$ on the right hand side of (\ref{phi}) is obtained by straightforward differentiation of (\ref{branch-mech}). \qed

 Equation (\ref{spine}) suggests that, under the
measure $\mathbb{P}^{-\lambda}_{\mu}$, the superprocess $X$ can be decomposed into two
parts.  The first one is a copy of the original superprocess and the
second one can be related to an independent process of immigration. As we shall demonstrate next, the process of immigration is governed by a spine or immortal particle along which two independent Poisson point processes of immigration occur. We need first to introduce some more notation.

Associated to the laws $\{\mathbb{P}_{\delta_x}:
x\in\mathbb{R}\}$ are the  measures $\{\mathbb{N}_x:
x\in \mathbb{R}\}$, defined on the same measurable space,  which
satisfy
\begin{equation}
\mathbb{N}_x (1- e^{-\langle f, X_t \rangle}) = -\log \mathbb{E}_{\delta_x}(e^{-\langle  f, X_t\rangle}),
\label{DK}
\end{equation}
for all $f\in C_b^+(\mathbb{R})$ and $t\geq 0$. Such measures are formally defined and explored in detail in \cite{Dynkin-Kuznetsov}. See also \cite{LG99}. Note that from the definition (\ref{DK}) it follows that
\[
\mathbb{N}_x(\langle f, X_t \rangle) = \mathbb{E}_{\delta_x}(\langle f, X_t \rangle),
\]
whenever $f\in C_b^+(\mathbb{R})$.

The measures $\{\mathbb{N}_x:x\in\mathbb{R}\}$ will play a crucial role in the forthcoming analysis. Intuitively speaking, the branching property implies that $\mathbb{P}_{\delta_x}$ is an infinitely divisible measure on the path space of $X$ and  (\ref{DK})  is a `L\'evy-Khinchine' formula in which  $\mathbb{N}_x$ plays the role of its `L\'evy measure'.  In this sense, $\mathbb{N}_x$ can be considered as the `rate' at which superprocesses `with zero initial mass' contribute to a unit mass at position $x$. It is important to note that $\mathbb{N}_x$ is not a probability measure as such.

With this measure in hand, let us now proceed to the definition of a
measure-valued process of immigration, which we denote by $\Lambda = \{\Lambda_t : t\geq 0\}$. Fix $x\in\mathbb{R}$ and $\mu\in\mathcal{M}_F(\mathbb{R})$.

      \begin{itemize}
\item[i.] ({\bf Spine}) We take a copy of the process
                  $\xi = \{\xi_t : t\geq 0\}$ under
                  $\Pi_x^{-\lambda}$ and henceforth refer to it as the  {\it  spine}.

\item[ii.]({\bf Continuum immigration})  Suppose that $\mathbf{n}$ is
                  a Poisson point process such that, for $t\geq 0$, given the
                  spine $\xi$, ${\bf n}$ issues a superprocess $X^{{\bf n},t}$ at space-time position $(\xi_t, t)$ with rate $ {\rm d}t\times 2 \beta{\rm d}\mathbb{N}_{\xi_t}$.

\item[iii.] ({\bf Jump immigration}) Suppose that $\mathbf{m}$ is a
            Poisson point process  such that, independently of ${\bf n}$, given the spine $\xi$,  $\mathbf{m}$ issues a superprocess $X^{{\bf m}, t}$ at space-time point $(\xi_t, t)$ with initial mass $r$ at rate
            ${\rm d}t\times r \nu({\rm d}r)\times {\rm d} \mathbb{P}_{r\delta_{\xi_t}}$. 


\end{itemize}
We now define for $t\geq 0$,
\begin{equation}\label{Lambda}
\Lambda_t = {X}'_t  +  X_t^{(\mathbf{n})} +  X_t^{(\mathbf{m})},
\end{equation}
where $\{X'_t: t\geq 0\}$ is an independent copy of $(X, \mathbb{P}_{\mu})$,
\begin{equation*}
 X_t^{(\mathbf{n})}=\sum_{s\leq t:\mathbf n}X^{\mathbf n,
 s}_{t-s},\, t\geq 0 \qquad
\mbox{and}\qquad
 X_t^{(\mathbf{m})}=\sum_{s\leq t:\mathbf{m}}X^{\mathbf m, s}_{t-s},\, t\geq 0.
\end{equation*}
In the last two equalities we understand the first sum to be over times at which the process ${\bf n}$ has a point and the second sum is understood similarly. Note that since the processes $ X^{(\mathbf{n})}$ and $ X^{(\mathbf{m})}$ are initially zero valued it is clear that since $X'_0 = \mu$ then $\Lambda_0 = \mu$. In that case, we use the notation $\widetilde{\mathbb{P}}^{-\lambda}_{\mu,x}$ to denote the law of the pair $( \Lambda, \xi)$. Note also that the pair $(\Lambda, \xi)$ are a time-homogenous Markov process. We are interested in the case that the initial position of the spine $\xi$ is randomised using the measure $\mu$ via (\ref{randomization}). In that case we shall write
\[
\widetilde{\mathbb{P}}^{-\lambda}_\mu(\cdot) = \frac{1}{\l e^{-\lambda \cdot}, \mu\r }\int_{\mathbb{R}} e^{-\lambda x}\mu({\rm d}x)\widetilde{\mathbb{P}}^{-\lambda}_{\mu, x}(\cdot)
\]
 for short. The next theorem identifies the process $\Lambda$ as the {\it pathwise} spine decomposition of $(X,\mathbb{P}^{-\lambda}_{\mu})$ and in particular it shows that as a process on its own $\Lambda$ is Markovian.

\begin{teo}\label{spine-decomp} For all $\lambda\in\mathbb{R}$ and $\mu\in\mathcal{M}_F(\mathbb{R})$, $(X, \mathbb{P}^{-\lambda}_{\mu})$ and $(\Lambda, \widetilde{\mathbb{P}}^{-\lambda}_{\mu})$ are equal in law.
\end{teo}
\proof Fix $\lambda\in\mathbb{R}$.
Firstly we must prove that
for any $t\geq 0$ and $\mu\in\mathcal{M}_F(\R)$,  we have
\begin{equation}
\mathbb{E}_\mu^{-\lambda}(e^{-\l g, X_t\r}) = \widetilde{\mathbb{E}}_\mu^{-\lambda}(e^{-\l g, \Lambda_t\r}),
\label{firstcheck}
\end{equation}
where $g\in C_b^+(\mathbb{R})$,
and for this it suffices to show that
\begin{equation}
\widetilde{\mathbb{E}}^{-\lambda}_\mu[e^{-\langle g,X_t^{(\mathbf{n})}+X_t^{(\mathbf{m})} \rangle}]
= \Pi_\mu^{-\lambda}\left[\exp\left\{-\int_0^t\phi(u_g(\xi_{t-s},s)){\rm d}s\right\}\right].
\label{firstthingtoprove}
\end{equation}
Secondly we must show that $(\Lambda, \widetilde{\mathbb{P}}^{-\lambda}_\mu)$ is a Markov process.

To this end note that for $g\in C_b^+(\mathbb{R})$,
\begin{eqnarray}\label{spine-one}
\lefteqn{ \widetilde{\mathbb{E}}^{-\lambda}_{\mu}\left[e^{-\l
 g,X^{(\mathbf{n})}_t+X^{(\mathbf{m})}_t\r}\right]}&&\notag\\
 &=&\nonumber\widetilde{\mathbb{E}}^{-\lambda}_{\mu}\left\{\widetilde{\mathbb{E}}^{-\lambda}_{\mu}\left[e^{-\l g,X^{(\mathbf{n})}_t\r}
                                     e^{-\l g,X^{(\mathbf{m})}_t\r}\bigg| \xi\right]\right\}\\
 &=&\Pi_\mu^{-\lambda}
              \bigg\{\widetilde{\mathbb{E}}^{-\lambda}_{\mu}\bigg[\exp\bigg\{-\sum_{s\leq
              t:\mathbf{n}}\l g,X^{\mathbf n,s}_{t-s}\r\bigg\}
              \bigg| \xi\bigg]
    \widetilde{\mathbb{E}}^{-\lambda}_{\mu}\bigg[\exp\bigg\{-
                              \sum_{s\leq t:\mathbf{m}}\l g,X^{\mathbf m, s}_{t-s}\r\bigg\}
              \bigg| \xi\bigg]\bigg\},
\end{eqnarray}
where we have used the independence of $X^{(\mathbf m)}$ and $X^{(\mathbf n)}$.
Applying Campbell's formula to the first inner expectation and
using \eqref{DK}, we get
\begin{eqnarray}\label{inner-one}
 \widetilde{\mathbb{E}}^{-\lambda}_{\mu}\left[\exp\{-\sum_{s\leq t:\mathbf{n}}
                                                \l g,X_{t-s}^{\mathbf n,s}\r\}\bigg| \xi\right]
 &=&\nonumber\exp\left\{-2\beta\int_0^t\int\left(1-e^{-\l g,X_{t-s}\r}\right){\rm d}\mathbb{N}_{\xi_s}{\rm d}s\right\}\\
 &=&\nonumber\exp\left\{-2\beta\int_0^t-\log \mathbb{E}_{\delta_{\xi_s}}\left(e^{-\l g,X_{t-s}\r}\right) {\rm d}s\right\}\\
 &=&\nonumber\exp\left\{-2\beta\int_0^tu_g(\xi_s,t-s){\rm d}s\right\}\\
 &=&\exp\left\{-2\beta\int_0^tu_g(\xi_{t-s},s){\rm d}s\right\}.
\end{eqnarray}
To deal with the second  expectation  in (\ref{spine-one}) first note that,
\begin{eqnarray}\label{inner-two}
 \widetilde{\mathbb{E}}^{-\lambda}_{\mu}\left[\exp\left\{-\sum_{s\leq t:\mathbf{m}}\l g,X^{\mathbf m,s}_{t-s}\r\right\}
            \bigg| \xi\right]
 = \widetilde{\mathbb{E}}^{-\lambda}_{\mu}\left[\exp\left\{-\sum_{s\leq t:\mathbf{m}}m_su_g(\xi_s,t-s)\right\}\bigg| \xi\right],
\end{eqnarray}
where for $s\geq 0$, $m_s = ||X^{\mathbf{m}, s}_0||$. In particular note that the process $\{m_t : t\geq 0\}$ is a Poisson point process on $(0,\infty)^2$, independent of $\xi$, with intensity ${\rm d} t \times r\nu({\rm d} r)$.
Hence, putting (\ref{inner-one}) and (\ref{inner-two}) into
(\ref{spine-one}) and again appealing to Cambell's formula yields
\begin{eqnarray*}
 \lefteqn{\widetilde{\mathbb{E}}^{-\lambda}_{\mu}\left[\exp\left\{-\l g,X^{(\mathbf{n})}_t+X^{(\mathbf{m})}_t\r\right\}\right]}&&\\
&=&\Pi_\mu^{-\lambda}\left\{\exp\left\{-2\beta\int_0^tu_g(\xi_{t-s},s){\rm d}s\right\}\widetilde{\mathbb{E}}^{-\lambda}_{\mu}\left[
             \exp\left\{-\sum_{s\leq t:{\mathbf{m}}}m_su_g(\xi_s,t-s)\right\}\bigg|\xi\right]\right\}\\
&=&\Pi_\mu^{-\lambda}\left\{\exp\left\{-2\beta\int_0^tu_g(\xi_{t-s},s){\rm d}s\right\}
                     \exp\left\{-\int_0^t\int_{(0,\infty)}\left(1-e^{-ru_g(\xi_{t-s},s)}\right)r\nu({\rm d}r){\rm d}s\right\}\right\}.
\end{eqnarray*}
Taking note of (\ref{phi}) we have in conclusion that
(\ref{firstthingtoprove}) follows.

Next we turn our attention to showing that $(\Lambda, \widetilde{\mathbb{P}}^{-\lambda}_\mu)$ is a Markov process. To this end, it suffices to show that for all $x\in\mathbb{R}$,
\begin{equation}
\widetilde{\mathbb{P}}^{-\lambda}_{\mu}(\xi_t\in {\rm d}x | \Lambda_t ) = \frac{1}{\langle e^{-\lambda\cdot}, \Lambda_t\rangle}e^{-\lambda x}\Lambda_t({\rm d}x).
\label{markovianproperty}
\end{equation}
Indeed, in that case it follows that for all $g\in C_b^+(\mathbb{R})$
\begin{eqnarray*}
\widetilde{\mathbb{E}}^{-\lambda}_{\mu}[e^{-\l g, \Lambda_t \r}| \Lambda_r : r\leq s] &=&
\widetilde{\mathbb{E}}^{-\lambda}_{\mu}\left[\left.\widetilde{\mathbb{E}}^{-\lambda}_{\mu}[e^{-\l g, \Lambda_t \r}| (\xi_r, \Lambda_r) : r\leq s]\right| \Lambda_r : r\leq s \right]\\
&=& \widetilde{\mathbb{E}}^{-\lambda}_{\mu}\left[\left. \left.\widetilde{\mathbb{E}}^{-\lambda}_{(\mu', x')}[e^{-\l g, \Lambda_t \r}]\right|_{ \mu' = \Lambda_s, x' = \xi_s}
\right| \Lambda_r : r\leq s \right]\\
&=& \widetilde{\mathbb{E}}^{-\lambda}_{\mu}\left[\left. \left.\widetilde{\mathbb{E}}^{-\lambda}_{( \mu', x')}[e^{-\l g, \Lambda_t \r}]\right|_{ \mu' = \Lambda_s, x' = \xi_s}
\right| \Lambda_s \right]\\
&=& \frac{1}{\l e^{-\lambda \cdot}, \Lambda_s\r}\int_{\mathbb{R}}e^{-\lambda x}\Lambda_s({\rm d}x)\left.\widetilde{\mathbb{E}}^{-\lambda}_{ \mu', x}[e^{-\l g, \Lambda_{t-s} \r}]\right|_{\mu' = \Lambda_s}\\
&=&\left. \widetilde{\mathbb{E}}^{-\lambda}_{\mu'}[e^{-\l g, \Lambda_{t-s} \r}]\right|_{\mu' = \Lambda_s},
\end{eqnarray*}
where in the third equality we have used that $(\xi,\Lambda)$ is Markovian.

To show (\ref{markovianproperty}) it suffices to show that for all $\theta\in \mathbb{R}$
\begin{equation}
\widetilde{\mathbb{E}}^{-\lambda}_{\mu}\left[e^{-\l g, \Lambda_t \r}
\widetilde{\mathbb{E}}^{-\lambda}_\mu[e^{-\theta \xi_t} |\Lambda_t]
\right]
= \widetilde{\mathbb{E}}^{-\lambda}_{\mu}\left[e^{-\l g, \Lambda_t \r}
\frac{\l e^{-(\lambda + \theta)\cdot}, \Lambda_t \r}{\l  e^{-\lambda \cdot} , \Lambda_t \r}
\right].
\label{sufftocheck}
\end{equation}
Note however that the left hand side of (\ref{sufftocheck}) is equal to
\begin{eqnarray}
\widetilde{\mathbb{E}}^{-\lambda}_{\mu}\left[
e^{-\l g, \Lambda_t \r}
e^{-\theta \xi_t} \right]
&=&e^{\gamma t}\frac{1}{\l e^{-\lambda \cdot}, \mu \r}\int_{\mathbb{R}} \mu({\rm d}x)e^{-(\lambda + \theta)x}\Pi_x^{-\lambda}\left[e^{- \theta(\xi_t - x) -\gamma t} \widetilde{\mathbb{E}}^{-\lambda}_\mu[e^{-\l g, \Lambda_t\r} | \xi_t]\right]\notag\\
&=&e^{\gamma t}\frac{1}{\l e^{-\lambda \cdot}, \mu \r}\int_{\mathbb{R}} \mu({\rm d}x)e^{-(\lambda + \theta)x}\Pi_x^{-(\lambda+\theta)}\left[ \widetilde{\mathbb{E}}^{-(\lambda+\theta)}_\mu[e^{-\l g, \Lambda_t\r} | \xi_t]\right]\notag\\
&=&e^{\gamma t}\frac{\l e^{-(\lambda+ \theta)\cdot}, \mu\r}{\l e^{-\lambda \cdot}, \mu \r}\widetilde{\mathbb{E}}^{-(\lambda+\theta)}_\mu[e^{-\l g, \Lambda_t\r} ],\label{1}
\end{eqnarray}
where $\gamma = (\lambda +\theta)^2/2 - \lambda^2/2$.

On the other hand, by (\ref{firstcheck}), the right hand side of (\ref{sufftocheck}) is equal to
\begin{eqnarray}
\mathbb{E}_{\mu}\left[e^{-\lambda c_\lambda t}\frac{\l e^{-\lambda \cdot}, X_t \r}{\l e^{-\lambda\cdot}, \mu \r}e^{-\l g, X_t \r}
\frac{\l e^{-(\lambda + \theta)\cdot}, X_t \r}{\l  e^{-\lambda \cdot} , X_t \r}
\right]&=&e^{\gamma t}\mathbb{E}_{\mu}\left[e^{-(\lambda +\theta)c_{\lambda + \theta}t}\frac{\l e^{-(\lambda+\theta) \cdot}, X_t \r}{\l e^{-\lambda\cdot}, \mu \r}e^{-\l g, X_t \r}\right]\notag\\
&=&e^{\gamma t}\frac{\l e^{-(\lambda+ \theta)\cdot}, \mu\r}{\l e^{-\lambda \cdot}, \mu \r}\mathbb{E}^{-(\lambda + \theta)}_{\mu}\left[e^{-\l g, X_t \r}\right].
\label{2}
\end{eqnarray}
Again appealing to (\ref{firstcheck}) we see that both (\ref{1}) and (\ref{2}) agree and the proof is complete.
\qed

\section{Proof of Theorem \ref{mgcgce} (i)}

For reasons of symmetry, it sufficies to prove the result with $\lambda\geq 0$.
Now that we are in possession of a pathwise spine decomposition, we may pursue a classical approach due to Lyons \cite{Ly}, see also Kyprianou \cite{Kyprianou2004},  to prove Theorem \ref{mgcgce} (i). The key element to the reasoning is the following measure theoretic result (see for example p242 of Durrett \cite{Durrett}).

Let
$\overline{W}_\infty(\lambda):=\limsup_{t\uparrow\infty}W_t(\lambda)$.
Then
\begin{equation}\label{Durret1}
 \overline{W}_\infty( \lambda)=\infty,\,\, \mathbb{P}^{-\lambda}\mbox{-a.s.}
 \Longleftrightarrow \overline{W}_\infty(\lambda)=0,\,\,\mathbb{P}\mbox{-a.s.}
\end{equation}
\begin{equation}\label{Durret2}
 \overline{W}_\infty( \lambda)<\infty,\,\, \mathbb{P}^{-\lambda}\mbox{-a.s.}
 \Longleftrightarrow \mathbb{E}(\overline{W}_\infty(\lambda))=1.
\end{equation}

Write for convenience $\widetilde{\mathbb{P}}^{-\lambda}$ instead of $\widetilde{\mathbb{P}}^{-\lambda}_{\delta_0}$.
Thanks to the spine decomposition (Theorem \ref{spine-decomp}) we may replace on the left hand sides of (\ref{Durret1}) and (\ref{Durret2}) $\mathbb{P}^{-\lambda}$  by $\widetilde{\mathbb{P}}^{-\lambda}$ and $\overline{W}_\infty(\lambda)$ by $\overline{W}^\Lambda_\infty (\lambda): =\limsup_{t\uparrow\infty} W^\Lambda_t (\lambda)$ where
\[
W^\Lambda_t (\lambda) :=  e^{-\lambda c_\lambda t}\langle e^{-\lambda \cdot}, \Lambda_t \rangle.
\]
We shall study $\overline{W}^\Lambda_\infty (\lambda) $ with the help of the martingale decomposition (which follows from the spine decomposition)
under $\widetilde{\mathbb{P}}^{-\lambda}$,
\begin{equation}\label{martingaledecompositon}
{W}^\Lambda_t (\lambda)=W'_t(\lambda)
 +\sum_{s\leq t:\mathbf{n}}e^{-\lambda c_\lambda s}W^{\mathbf n,s}_{t-s}(\lambda)
 +\sum_{s\leq t:\mathbf{m}}e^{-\lambda c_\lambda s}W^{\mathbf m,s}_{t-s}(\lambda),
\end{equation}
where $W'(\lambda)$ is an independent copy of $W(\lambda)$ under $\mathbb{P}$,  $W^{\mathbf{n},s}_{t-s}(\lambda) = e^{ - \lambda c_\lambda(t-s)}\langle e^{-\lambda \cdot}, X^{\mathbf{n}, s}_{t-s}\rangle$
and $W^{\mathbf{m},s}_{t-s}(\lambda) = e^{ - \lambda c_\lambda (t-s)}\langle e^{-\lambda \cdot}, X^{\mathbf{m}, s}_{t-s}\rangle$. 

\bigskip

Suppose  that $\int_{[1,\infty)} r(\log r)\nu({\rm d}r)<\infty$ and
$\lambda\in(0,\underline{\lambda})$. First note that
since $W(\lambda)$ is a $\mathbb{P}$-martingale, it follows as a
standard result that $W(\lambda)^{-1}$ is a
$\mathbb{P}^{-\lambda}$-supermartingale (see for example \cite{HR})
and hence $\lim_{t\uparrow\infty}W_t(\lambda)$ exists
$\mathbb{P}^{-\lambda}$-almost surely, and hence
$\lim_{t\uparrow\infty} W^\Lambda_t (\lambda)$ exists
$\widetilde{\mathbb{P}}^{-\lambda}$-almost surely and is equal to
$W^\Lambda_\infty (\lambda) $. Our objective is to show that
\begin{equation}
\limsup_{t\uparrow\infty}\widetilde{\mathbb{E}}\left[\left. W^\Lambda_t (\lambda) \right|\xi, {\bf m}\right]<\infty,
\label{objective}
\end{equation}
in which case Fatou's lemma and the existence of  $\lim_{t\uparrow\infty} W^\Lambda_t (\lambda) = \overline{W}^\Lambda_\infty (\lambda) $ implies that, $\widetilde{\mathbb{P}}^{-\lambda}$-almost surely, $\overline{W}^\Lambda_\infty (\lambda) <\infty$. This in turn implies that $\mathbb{P}^{-\lambda}$-almost surely, $\overline{W}_\infty(\lambda) <\infty$ and hence, by (\ref{Durret2}), $W_\infty(\lambda)$ is an $L^1(\mathbb{P})$ limit.



To this end note that, since $\widetilde{\mathbb{E}}^{-\lambda}[W'_t(\lambda)] = 1$, it suffices to prove that
\[
\limsup_{t\uparrow\infty}
 \widetilde{\mathbb{E}}^{-\lambda}\left[\left.\sum_{s\leq t:\mathbf{n}}
 e^{-\lambda c_\lambda s}W_{t-s}^{\mathbf n, s}(\lambda) + \sum_{s\leq t:\mathbf{m}}
 e^{-\lambda c_\lambda s}W_{t-s}^{\mathbf m, s}(\lambda)\right|\xi, {\bf m}\right]<\infty.
\]
First note that, $\widetilde{\mathbb{P}}^{-\lambda}$-almost surely,
\begin{eqnarray}\label{cont-part}
\limsup_{t\uparrow\infty}    \widetilde{\mathbb E}^{-\lambda}\left[
    \sum_{s\leq t:\mathbf{n}}e^{-\lambda c_\lambda s}W^{\mathbf n, s}_{t-s}(\lambda)\bigg|\xi,\mathbf{m}\right]
 &=&\limsup_{t\uparrow\infty}\nonumber 2\beta \int_0^te^{-\lambda c_\lambda s}\mathbb{N}_{\xi_s}[W_{t-s}(\lambda)] {\rm d}s\\
 &=&\nonumber\limsup_{t\uparrow\infty} 2\beta \int_0^te^{-\lambda c_\lambda s}\mathbb{E}_{\delta_{\xi_s}}
    \left[W_{t-s}(\lambda)\right]{\rm d}s\\
 &=&2\beta \int_0^\infty e^{-\lambda(\xi_s+c_\lambda s)}{\rm d}s\notag\\
 &<&\infty,
\end{eqnarray}
where the second identity holds by (\ref{DK}) and the final inequality is a result of the strong law of large numbers for linear Brownian motion together with the fact that $\xi_t + c_\lambda t$ has
drift $c_\lambda-\lambda$ which is strictly positive when $\lambda\in(0,\underline{\lambda})$. 

Next, recalling that for all $x\in\mathbb{R}, t\geq 0, r>0$, $\mathbb{E}_{r\delta_x}(W_t(\lambda)) = r e^{-\lambda x}\mathbb{E}(W_t(\lambda)) =  r e^{-\lambda x}$, we have that
\begin{eqnarray}\label{cod-exp1}
 \widetilde{\mathbb{E}}^{-\lambda}\left[\left.\sum_{s\leq t:\mathbf{m}}
 e^{-\lambda c_\lambda s}W_{t-s}^{\mathbf m, s}(\lambda)\right|\xi, \mathbf{m}\right]
 &=&\widetilde{\mathbb{E}}^{-\lambda}\left[\left.\sum_{s\leq t:\mathbf{m}}m_se^{-\lambda(\xi_s+c_\lambda
 s)}\right|\xi,\mathbf{m}\right]\notag\\
& \leq&
\sum_{s\geq 0;\mathbf{m}}m_se^{-\lambda(\xi_s+c_\lambda
 s)},
\end{eqnarray}
where we recall that the process $\{m_t: t\geq 0\}$ is a Poisson point process, independent of $\xi$ and with intensity ${\rm d}t\times r\nu({\rm d}r)$
However, on the one hand we have that $\widetilde{\mathbb{P}}^{-\lambda}$-almost surely
\begin{eqnarray}\label{finitesmall1}
 \widetilde{\mathbb E}^{-\lambda}
     \left[\sum_{s\geq 0;\mathbf{m}}m_{s}1_{\{m_s<1\}}e^{-\lambda(\xi_s+c_\lambda s)}\bigg|\xi
     \right]
 &=&\nonumber\int_{(0,1)}\int_0^\infty re^{-\lambda(\xi_s+c_\lambda s)}r{\rm d}s\nu({\rm d}r)\\
 &=&\nonumber\int_{(0,1)}r^2d\nu({\rm d}r)\int_0^\infty e^{-\lambda(\xi_s+c_\lambda s)}{\rm d}s\\
 &<&\infty.
\end{eqnarray}
On the other hand,  define $\tau_0=0$ and $\tau_i=\inf\{t\geq\tau_{i-1}: m_t\geq1\}$,
$i=1,2,\cdots$.  Note that, $\{\tau_i: i\geq 0\}$ are the times or arrival in  a Poisson process with
arrival rate $\int_{[1,\infty)} r\nu({\rm d}r)$.
Recall that  if $\{Z_i: i\geq 1\}$ is a
sequence of i.i.d. random variables with probability measure $Q$
then it can easily be shown with the help of the Borel-Cantelli Lemma that, $Q$-a.s.
\begin{equation}\label{result}
\limsup_{i\rightarrow\infty} i^{-1}\log Z_i = \left\{
\begin{array}{ll}
         0, & \mbox{if}\,\, Q(\log^+Z_1)<\infty;\\
       \infty, & \mbox{if}\,\,Q(\log^+Z_1)=\infty.
       \end{array} \right.
\end{equation}
In this instance we would like to take $Z_i = m_{\tau_i}$ which has common probability
measure
$Q({\rm d}r)=r\nu({\rm d}r)/\int_{[1,\infty)} r\nu({\rm d}r)
$
on $[1,\infty)$.
We have  $\widetilde{\mathbb{P}}^{-\lambda}$-almost surely that
\begin{equation}\label{bigm}
\sum_{s\geq 0;\mathbf{m}}m_{s}1_{\{m_s\geq 1\}}e^{-\lambda(\xi_s+c_\lambda s)}
 =\sum_{i\geq1}m_{\tau_i}e^{-\lambda(\xi_{\tau_i}+c_\lambda \tau_i)}\\
 <\infty,
\end{equation}
where the equality follows by virtue of the fact that entries in the first sum arrive at rate $\int_{[1,\infty)}r\nu({\rm d}r)$, which is finite,  and the final inequality follows from (\ref{result}), the assumption $\int_{[1,\infty)} r\log r\nu({\rm d}r)<\infty$, the
 fact that
\begin{equation}
\lim_{i\uparrow\infty}\frac{\tau_i}{i} = \frac{1}{\int_{[1,\infty)} r\nu({\rm d}r)},
\label{SLLNtau}
\end{equation}
and the strong law of large numbers for linear Brownian motion.
In conclusion, (\ref{finitesmall1}) and (\ref{bigm}) show that (\ref{cod-exp1}) is finite and hence $W_\infty(\lambda)$ is an $L^1(\mathbb{P})$-limit.

 Finally we prove that $\{W_\infty(\lambda) = 0\}$ agrees with the event that $X$ becomes extinguished. To this end let $q= \mathbb{P}_{\delta_x}(W_\infty(\lambda) = 0)$. Note that $q$ does not depend on $x$ as $W_\infty(\lambda)$ under $\mathbb{P}_{\delta_x}$ has the same law as $e^{-\lambda x}W_{\infty}(\lambda)$ under $\mathbb{P}_{\delta_0}$ thanks to the definition of $W(\lambda)$ and the fact that the branching mechanism is not spatially dependent.
Taking conditional expectations and using the Markov branching property we have for all $t\geq 0$,
\[
\mathbb{E}(\mathbf{1}_{\{W_\infty(\lambda) = 0\}}| \mathcal{F}_t)
= q^{||X_t||}.
\]
Therefore, on the one hand, taking expectations across this last equality and then limits as $t\uparrow\infty$, we easily deduce with the help of the Dominated Convergence Theorem that $\mathbb{P}(W_\infty(\lambda) = 0) = e^{-\lambda^*}$. On the other hand, noting that $\mathcal{E}^c\in\sigma\left(\bigcup_{t>0}\mathcal{F}_t\right)$,  we also have that
\[
\mathbb{P}(\mathcal{E}^c\cap \{W_\infty (\lambda) = 0\}) = \mathbb{E}\left[
\mathbf{1}_{\mathcal{E}^c}\lim_{t\uparrow\infty}\mathbb{E}(\mathbf{1}_{\{W_\infty(\lambda) = 0\}}| \mathcal{F}_t)
\right] =\mathbb{E}[\mathbf{1}_{\mathcal{E}^c} \lim_{t\uparrow\infty}q^{||X_t||}] =0.
\]
Hence it follows that $\{W_\infty (\lambda) = 0\} = \mathcal{E}$, $\mathbb{P}$-almost surely.


\bigskip

Next we deal with the cases that
$\lambda\geq \underline{\lambda}$ or $\int_{[1,\infty)} r(\log r)\nu({\rm d}r)=\infty$. Recall that, given $\xi$, the Poisson point process $\mathbf{m}$ initiates a superprocess at time $t$,  $X^{\mathbf{m}, t}$, with  initial mass $r$ at rate ${\rm d}t\times r\nu({\rm d}r)\times {\rm d}\mathbb{P}_{r\delta_{\xi_t}}$.
 For each $\tau_i$ (defined in the previous part of the proof) we have
\begin{equation}\label{limit-martingale-1}
W^\Lambda_{\tau_i} (\lambda) \geq
m_{\tau_i}e^{-\lambda(\xi_{\tau_i}+c_\lambda\tau_i)},\qquad\widetilde{\mathbb{P}}^{-\lambda}-\mbox{a.s.}
\end{equation}
Under $\widetilde{\mathbb P}$, $\xi$ is a
Brownian motion with drift $-\lambda$ and is independent of
$\mathbf{m}$.  Note also  that $c_\lambda\leq\lambda$ when
$\lambda\geq\underline{\lambda}$.  Hence, $\xi_t+c_\lambda
t$ is a Brownian motion with non-positive drift.     Then, from
\eqref{limit-martingale-1} and (\ref{SLLNtau}) we conclude
that when $\lambda\geq \underline{\lambda}$,
\begin{equation*}
     \limsup_{t\rightarrow\infty} W^\Lambda_t(\lambda)
 \geq\limsup_{i\rightarrow\infty} e^{-\lambda(\xi_{\tau_i}+c_\lambda \tau_i)}
    =\infty,\,\,\widetilde{\mathbb P}-\mbox{a.s.}
\end{equation*}
It follows from (\ref{Durret1}) that $W_\infty(\lambda)  = 0$, $\mathbb{P}$-almost surely.

\bigskip

Now suppose that $\int_{[1,\infty)} r(\log r)\nu({\rm d}r)=\infty$.
Recall that under $\widetilde{\mathbb{P}}^{-\lambda}$, $\xi$ is a Brownian
motion with drift $- \lambda$.  Then, the strong law of large
numbers gives  $\widetilde{\mathbb{P}}^{-\lambda}$-almost surely,
 $
\lim_{t\rightarrow\infty}t^{-1}(\xi_t+ c_\lambda t)=c_\lambda-\lambda.
 $
  We thus have from (\ref{result}), (\ref{SLLNtau}) and (\ref{limit-martingale-1}) that
$\widetilde{\mathbb{P}}^{-\lambda}$-a.s.,
\begin{equation}
     \limsup_{i\rightarrow\infty}i^{-1}{\log W^\Lambda_{\tau_i}(\lambda)}
 \geq\limsup_{i\rightarrow\infty}i^{-1}{\log m_{\tau_i}}
    -\lambda\limsup_{i\rightarrow\infty}\frac{(\xi_{\tau_i}+c_\lambda\tau_i)
     }{\tau_i}\frac{\tau_i}{i}
    =\infty.
\end{equation}
Therefore, $\overline{W}^\Lambda_\infty(\lambda)=\infty$,
$\widetilde{\mathbb{P}}^{-\lambda}$-almost surely, and (\ref{Durret1}) implies that
$W_\infty(\lambda)=\overline{W}_\infty(\lambda)=0$, $\mathbb{P}$-almost surely.
\qed

\begin{remark}\rm
One can easily go further with the analysis of $W(\lambda)$ through the pathwise spine decomposition. Indeed following arguments in Hardy and Harris \cite{HH},  sufficient conditions for $L^p(\mathbb{P})$ convergence for $p\in(1,2]$ are given in \cite{Lp}.
\end{remark}

\section{Proof of Theorem \ref{mgcgce} (ii)}

As with the case of $W(\lambda)$ we only give the proof for the case that $\lambda\geq 0$ noting that
reasoning involving symmetry covers the case that $\lambda\leq 0$.
Recall that $\partial W(\lambda)$ is a signed martingale and therefore does not necessarily converge
almost surely. A technique used by Kyprianou \cite{Kyprianou2004} to get round this problem in the
case of a branching Brownian motion is to consider a truncated form of the derivative martingale which is a positive martingale.
In order to describe the aforementioned martingale in the current context we need more notation.

Recall that the superprocess $X^{c}$ was defined as
the superprocess whose movement component is that of a Brownian
motion with drift $c$ but whose branching mechanism is
still $\psi$. Consider the domain $D^{t}_{-y}  = \{(z,s): z>-y, \,
s\in(0,t)\}$. Dynkin's theory of exit measures may still be applied
in this context and we denote the exit measure associated with the
domain $D^{t}_{-y}  $ for the process $X^{\lambda}$ by
$X^{\lambda}_{D^{t}_{-y}  }$ for $\lambda\geq 0$. Next define for
all $\lambda\geq 0$,
\[
b_\lambda  = c_\lambda-\lambda = -\psi'(0^+)/\lambda-\lambda/2,
\]
and note that $b_\lambda>0$ for $\lambda\in
(0,\underline\lambda)$ and $b_\lambda\leq 0$ for
$\lambda\geq\underline\lambda$.

In the spirit of \cite{Kyprianou2004} we introduce a new  martingale,  for each $y>0$,
\begin{eqnarray}\label{second martingale}
V^{-y}_t(\lambda)= e^{-\lambda b_\lambda t}\frac{1}{y}\langle  (y+\cdot)e^{-\lambda\cdot}, X^{\lambda}_{D^{t}_{-y}  } \rangle,\,\, t\geq 0.
\end{eqnarray}
To show that $V^{-y}(\lambda) : = \{V^{-y}_t(\lambda) : t\geq 0\}$
is a $\mathbb{P}$-unit mean martingale, let
$\mathcal{H}^{\lambda}_t =
\sigma(X^{\lambda}_{D^{s}_{-x}} : x\leq y,  s\leq t )\subseteq {\color{black}{\mathcal{K}}_t=\sigma(X^{\lambda}_{D^{s}_{-x}} : x<\infty,  s\leq t )}$
and note  that
\[
\mathbb{E}(y W_t(\lambda) + \partial W_t(\lambda)| \mathcal{H}^{\lambda}_t) =y V^{ -y}_t(\lambda),
\]
and that $\mathcal{H}^{\lambda}_{t+s}\cap  {\color{black}{\mathcal{K}}_t}  = \mathcal{H}^{\lambda}_t$.
Hence
\begin{eqnarray*}
\mathbb{E}(yV^{ -y}_{t+s} (\lambda)| {\color{black}{\mathcal{K}}_t}) &=&
 \mathbb{E}(y W_{t+s}(\lambda) + \partial W_{t+s}(\lambda)| \mathcal{H}^{\lambda}_{t+s}| {\color{black}{\mathcal{K}}_t})\\
 &=& \mathbb{E}(y W_{t+s}(\lambda) + \partial W_{t+s}(\lambda) |\mathcal{H}^{\lambda}_t)\\
 &=& \mathbb{E}(y W_{t+s}(\lambda) + \partial W_{t+s}(\lambda)| {\color{black}{\mathcal{K}}_t} |\mathcal{H}^{\lambda}_t)\\
 &=&\mathbb{E}(y W_{t}(\lambda) + \partial W_{t}(\lambda)|\mathcal{H}^{\lambda}_t)\\
 &=&yV^{ -y}_t(\lambda).
\end{eqnarray*}

It is clear that $V^{-y}(\lambda)$ is  positive and hence there always exists an almost sure limit which we denote by $V^{-y}_\infty(\lambda)$. {\color{black} From} Corollary \ref{support} we know that, when $\lambda\geq \underline\lambda$,  $\mathbb{P}(-\inf\mathcal{R}^{\lambda}<\infty)=1$. It follows that on the event $\{\inf\mathcal{R}^{\lambda}\geq -y\}$, for this regime of $\lambda$,
\[
\mathbb{E}(y W_t(\lambda) + \partial W_t(\lambda)| \mathcal{H}^{\lambda}_t) =
y W_t(\lambda) + \partial W_t(\lambda),
\]
and hence, letting $t\uparrow\infty$,  $yV^{-y}_\infty(\lambda) = yW_\infty(\lambda) + \partial W_\infty(\lambda)$, where implicitly we understand the limit of $\partial W(\lambda)$ to exist in the last equality because the limit $V^{-y}_\infty(\lambda)$ exists. Note however from Theorem \ref{mgcgce} (i) that $W_\infty(\lambda)=0$ when $\lambda\geq \underline\lambda$ so that in fact
\begin{equation}
yV^{-y}_\infty(\lambda) =  \partial W_\infty(\lambda) \text{ on } \{\inf\mathcal{R}^{\lambda} \geq -y\}.
\label{V=W}
\end{equation}
 As $y$ may be taken arbitrarily large, it follows that $\partial W_\infty (\lambda)\geq 0$.

 Remaining in the regime $\lambda\geq \underline\lambda$, the proof that  $\{\partial W_\infty (\lambda) = 0\} = \mathcal{E}$, $\mathbb{P}$-almost surely goes along almost the same lines as the earlier proof that $\{W_\infty (\lambda) = 0\} = \mathcal{E}$, $\mathbb{P}$-almost surely.


\bigskip

Taking account of the relationship between $\partial W(\lambda)$ and $V^{-y}(\lambda)$ for $\lambda\geq \underline\lambda$, the proof of Theorem \ref{mgcgce} (ii) would now follow directly from parts (ii) and (iii) of the following theorem; which itself plays the analogous role of Theorem 13 in Kyprianou \cite{K1}.

\begin{theorem}\label{13}
Fix $y>0$.
\begin{itemize}
\item[(i)]  If $\lambda > \underline\lambda$ then $V^{-y}_\infty (\lambda) = 0$ $\mathbb{P}$-almost surely.
\item[(ii)] If $\lambda = \underline\lambda$ then $V^{-y}_\infty (\lambda)$ is an $L^1(\mathbb{P})$ limit if and only if $\int_{[1,\infty)}r (\log r)^2\nu({\rm d}r)<\infty$ otherwise $V^{-y}_\infty (\lambda) = 0$ $\mathbb{P}$-almost surely.
\item[(iii)] If $\lambda\in(0,\underline\lambda)$ then $V^{-y}_\infty (\lambda)$ is an $L^1(\mathbb{P})$ limit if and only if $\int_{[1,\infty)} r(\log r)\nu({\rm d}r)<\infty$ otherwise $V^{-y}_\infty (\lambda) = 0$ $\mathbb{P}$-almost surely.
\end{itemize}
\end{theorem}

Below we only give  the proof of part (ii) of Therorem \ref{13}. Once this has been done, the proof of parts (i) and (iii) should be apparent given the proof of Theorem 13 in Kyprianou \cite{K1} and we leave the details to the reader.

\bigskip

\noindent {\bf Proof of Theorem \ref{13} (ii):}
We shall again appeal to classical techniques based around using a martingale  change of measure.
Specifically we are interested in understanding the change of measure
\begin{equation}
\label{phat}
\left.\frac{{\rm d}\widehat{\mathbb{P}}^{-y} }{{\rm d}\mathbb{P}}\right|_{\mathcal{F}_{t}} := V^{-y}_t(\underline\lambda), \,\,\, t\geq 0
\end{equation}
where $y>0$. Similarly to Theorem \ref{spine-decomp} the change
of measure induces a spine decomposition. In order to describe it,
recall that under $\Pi_x$ the process $\xi:=\{\xi_t:t\geq 0\}$ is a
Brownian motion issued from $x\in\mathbb{R}$. If we let
 $\tau_{-y} = \inf\{t\geq 0 : y+\xi_t+\underline\lambda
t\leq 0\}$ then another well known change of measure for Brownian
motion is the following. For $y\geq 0$,
\[
\left.\frac{ {\rm d}\widehat\Pi^{-y}}{{\rm d}\Pi}\right|_{\mathcal{G}_t}
:= \frac{y+\xi_t +{\underline{\lambda}} t}{y}e^{-{\underline{\lambda}}\xi_t - {\underline{\lambda}}^2 t/2} \mathbf{1}_{\{t<\tau_{-y}\}}, \,\, t\geq 0
\]
where $\mathcal{G}_t = \sigma(\xi_s: s\leq t)$ and
$\Pi = \Pi_0$. Under $\widehat\Pi^{-y}$ the process $\{y+\xi_t
+\underline{\lambda} t: t\geq 0\}$ has the law of a standard
Brownian motion issued from $y$ and conditioned never to enter the
half line $(-\infty,0)$. Otherwise said, the process $y+\xi_t +
\underline{\lambda}t$ is  a Bessel-3 process issued from $y$.
Bearing this last change of measure in mind, we have the following
result which describes the effect of the change of measure
(\ref{phat}).

\begin{theorem}\label{second spine} Fix $y\geq 0$.
Consider the process $\Lambda$ as defined in (\ref{Lambda}) with the
exception that the spine $\xi$ is assigned the measure
$\widehat\Pi^{-y}$  such that
$\{y+\xi_t+\underline\lambda t: t\geq 0\}$ is a Bessel-3 process issued from $y$ and $x$ is chosen specifically equal to 0.
Denote its law by $\widetilde{\mathbb{P}}^{-y}$. Then $(X,
\widehat{\mathbb{P}}^{-y}) = (\Lambda,
\widetilde{\mathbb{P}}^{-y})$.
\end{theorem}

For the sake of brevity we omit the proof mentioning instead that it requires  very similar computations, with obvious differences, to those of Theorems \ref{spineth} and \ref{spine-decomp} combined.

\bigskip

Theorem \ref{second spine} allows us to conclude that the process
$V^{-y}(\lambda)$ under $\widehat{\mathbb{P}}^{-y}$ is equal in law
to
\begin{eqnarray}\label{cri-decom}
 V^{\Lambda, -y}_t(\underline\lambda) &: =& V^{\prime -y}_t(\underline\lambda)+\sum_{s\leq t:\mathbf{n}}\frac{(y + \xi_s + \underline{\lambda} s)}{y} e^{-\underline\lambda(\xi_s +\underline\lambda s)} V_{t-s}^{{\mathbf n}, s, -(y + \xi_s + \underline{\lambda} s)}(\underline\lambda)\nonumber\\
&& +\sum_{s\leq t:\mathbf{m}}  \frac{(y + \xi_s +
\underline{\lambda} s)}{y}e^{-\underline \lambda(\xi_s
+\underline\lambda s)}V_{t-s}^{{\mathbf m}, s, -(y + \xi_s +
\underline{\lambda} s)}(\underline\lambda), \hspace{1.5cm} t\geq
0,
\end{eqnarray}
under $\widetilde{\mathbb{P}}^{-y}$, where $V^{\prime,
-y}(\underline\lambda) $ plays the role of
$V^{-y}(\underline\lambda)$ for the process $X'$ and, given $(\xi,
{\bf m})$, $V^{{\mathbf m}, s, -(y + \xi_s + \underline{\lambda}
s)}(\underline\lambda)$ and  $V^{{\mathbf n}, s, -(y + \xi_s +
\underline{\lambda} s)}(\underline\lambda)$ play the role of $V^{-(y
+ \xi_s + \underline{\lambda} s)}(\underline\lambda)$ for the
processes $X^{{\mathbf m}, s}$  and $X^{{\mathbf n}, s}$,
respectively under the laws $\mathbb{P}_{m_s\delta_{0}}$ and
$\mathbb{N}_0$.

A similar statement  to  (\ref{Durret2}) tells us that if we can show $\widetilde{\mathbb{P}}^{-y}(\limsup_{t\uparrow\infty}  V^{\Lambda, -y}_t(\underline\lambda) <\infty) = 1$ then $V^{-y}_\infty(\underline\lambda)$ is an  $L^1({\mathbb{P}})$ limit.
Similar reasoning to  the  proof of Theorem \ref{mgcgce} (i) tells us that it now suffices to prove that
\begin{equation}
\limsup_{t\uparrow\infty} \widetilde{\mathbb{E}}^{-y}(V^{\Lambda, -y}_t (\underline\lambda)|\xi,{\bf m})<\infty,
\label{Vlimsup}
\end{equation}
almost surely. To this end, first recall that, given $\xi$, the Poisson point process $\mathbf{n}$ of immigration has intensity  $ {\rm d}s\times2\beta{\rm d}\mathbb{N}_{\xi_s}$.
It follows that
\begin{eqnarray*}
\lefteqn{
\limsup_{t\uparrow\infty}\widetilde{\mathbb{E}}^{-y}\left[\left.\sum_{s\leq t:\mathbf{n}} (y + \xi_s + \underline{\lambda} s) e^{-\underline\lambda(\xi_s +\underline\lambda s)}V_{t-s}^{{\mathbf n}, s, -(y + \xi_s + \underline{\lambda} s)}(\underline\lambda)
\right| \xi\right]}&&\\
&=&\limsup_{t\uparrow\infty} 2\beta \int_0^t  (y + \xi_s + \underline{\lambda} s)e^{-\underline\lambda(\xi_s +\underline\lambda s)} \mathbb{N}_0(V_{t-s}^{{\mathbf n}, s, -(y + \xi_s + \underline{\lambda} s)}(\underline\lambda)|\xi){\rm d}s \\
&=&\limsup_{t\uparrow\infty} 2\beta \int_0^t  (y + \xi_s + \underline{\lambda} s)e^{-\underline\lambda(\xi_s +\underline\lambda s)} \mathbb{E}( V_{t-s}^{ -(y + \xi_s + \underline{\lambda} s)}(\underline\lambda)|\xi){\rm d}s \\
&=&\int_0^\infty  (y + \xi_s + \underline{\lambda} s) e^{-\underline{\lambda}( \xi_s +\underline{\lambda} s) }{\rm d}s\\
&<&\infty,
\end{eqnarray*}
$\widetilde{\mathbb{P}}^{-y}$-almost surely,
where the final equality follows by virtue of the fact that $\{y+\xi_t + \underline{\lambda} t: t\geq 0\}$ is a Bessel-3 process issued from $y$ and hence eventually grows no slower than $t^{\frac{1}{2}-\epsilon}$ for any $1/2>\epsilon>0$.

Next we note that
\begin{eqnarray}
\lefteqn{\hspace{-3cm}
\widetilde{\mathbb{E}}^{-y} \left[
\left.
\sum_{s\leq t:\mathbf{m}}  (y + \xi_s + \underline{\lambda} s) e^{-\underline\lambda(\xi_s +\underline\lambda s)}V_{t-s}^{{\mathbf m}, s, -(y + \xi_s + \underline{\lambda} s)}(\underline\lambda)
\right|\xi,\mathbf{m}
\right]
}
&&\\
&&=
\sum_{s\leq t:\mathbf{m}}  (y + \xi_s + \underline{\lambda} s)e^{-\underline\lambda(\xi_s +\underline\lambda s)}\mathbb{E}_{m_s\delta_0}(V_{t-s}^{{\mathbf m}, s, -(y + \xi_s + \underline{\lambda} s)}(\underline\lambda)|\xi,\mathbf{m}) \notag\\
&&=\sum_{s\leq t:\mathbf{m}}m_s   (y + \xi_s + \underline{\lambda} s)e^{-\underline{\lambda}( \xi_s+\underline{\lambda} s)}\notag\\
&&= \sum_{\{m_s<e^{\varepsilon(\xi_s+{\underline{\lambda}}
   s)}:\mathbf m\}} m_s (y+\xi_s+{\underline{\lambda}} s) e^{-{\underline{\lambda}}(\xi_s+{\underline{\lambda}} s)}
\notag\\
&&+\sum_{\{m_s\ge
   e^{\varepsilon(\xi_s+{\underline{\lambda}}
   s)}:\mathbf m\}} m_s  (y+\xi_s+{\underline{\lambda}} s)e^{-{\underline{\lambda}}(\xi_s+{\underline{\lambda}} s)}\notag\\
&&=: I + II.
\label{I+II}
\end{eqnarray}
We want to show that $I$ and $II$ are
both  $\widetilde{\mathbb{P}}^{-y}$-almost surely
finite. For $I$, choose $0<\varepsilon<{\underline{\lambda}}$, we
have that
\begin{eqnarray*}
  \widetilde{\mathbb{E}}^{-y}(I |\xi) &\leq &\int_{(0,1)}r^2\nu({\rm d}r)\int_0^\infty (y+\xi_t+{\underline{\lambda}}
  t)   e^{-{\underline{\lambda}}(\xi_t+{\underline{\lambda}} t)}{\rm d}t
\\
  &+&
\int_{[1,\infty)} r\nu({\rm d}r)\int_0^\infty
(y+\xi_t+{\underline{\lambda}} t)e^{-({\underline{\lambda}}-\varepsilon)(\xi_t+{\underline{\lambda}} t)}
{\rm d}t\\
&<&\infty.
\end{eqnarray*}
Note that we have again used the fact that the assumption $\psi'(0^+)\in(-\infty, 0)$ implies that $\int_{[1,\infty)} r\nu({\rm d}r)<\infty$.

{\color{black}
  To show that $II$ is $\widetilde{\mathbb{P}}^{-y}$-almost surely
finite, it suffices to note that
\begin{eqnarray*}
\lefteqn{
\widetilde{\mathbb{E}}^{-y}\left[
\sum_{\{m_s\ge
   e^{\varepsilon(\xi_s+{\underline{\lambda}}
   s)}:\mathbf m\}} 1
\right]
}&&\\
&=& \widehat\Pi^{-y}_0\int_0^\infty {\rm d}t\int_{[1,\infty)}
r\nu({\rm d}r)1_{\{r>e^{\varepsilon(\xi_t+{\underline{\lambda}}
   t)}\}}\\
&=&\widehat\Pi^{-y}_0\int_{[1,\infty)} r\nu({\rm d}r)\int_0^\infty
{\rm d}t1_{\{x+\xi_t+{\underline{\lambda}} t\leq x+\varepsilon^{-1}\log
   r\}}\\
   &=&\int_{[1,\infty)} r\nu({\rm d}r)\int_{\{|{\bf y}|\leq x+\varepsilon^{-1}\log
  r\}}{\rm d}{\bf y}\int^\infty_0 p(t, \hat{{\bf x}},{\bf y}) {\rm d}t\\
&=&C\int_{[1,\infty)} r\nu({\rm d}r)\int_{\mathbb
R^3}\dfrac{1_{\{|{\bf y}|\leq x+\varepsilon^{-1}\log
  r\}}}{|{\bf y}-\hat{\bf x}|}{\rm d}{\bf y}\\
&\leq&C\int_{[1,\infty)}
r\nu({\rm d}r)\int_0^{2x+\varepsilon^{-1}\log
  r}u{\rm d}u\\
 &=&\frac{C}{2}\int_{[1,\infty)} r(2x+\varepsilon^{-1}\log r)^2\nu({\rm d}r)\\
 &<&\infty,
\end{eqnarray*}
where $\hat{{\bf x}}=(x,0,0)$, $p(t, \hat{{\bf x}},{\bf y})$ is the probability density function of a three dimensional Brownian motion starting from $\hat{{\bf x}}$, and $C$ is a positive constant. This tells us that, $\widetilde{\mathbb{P}}^{-y}$-almost surely,  $II$ is a summation over a finite set, and therefore $II$ is $\widetilde{\mathbb{P}}^{-y}$-almost surely
finite.
}


 A similar statement to  (\ref{Durret1}) tells us that if we can
show $\widetilde{\mathbb{P}}^{-y}(\limsup_{t\uparrow\infty}
V^{\Lambda, -y}_t(\underline\lambda) =\infty) = 1$ then
 $V^{-y}_\infty(\underline\lambda)=0$ ${\mathbb{P}}$- almost surely.  Since each term in the identity
  \eqref{cri-decom} is nonnegative, we just consider the $m$-immigration.  If we prove the supremum limit of the third term
  in \eqref{cri-decom} is infinity, then we are done.  Let $N$ be any positive
number and define the stochastic time sequence
$$
\tau_1=\inf\{t\geq 0: m_t>1\vee e^{N(\xi_t+\underline\lambda t)}\},
\tau_{i+1}=\inf\{t>\tau_i: m_t>1\vee e^{N(\xi_t+\underline\lambda
t)}\}, i=1,\ldots
$$
 Under $\widetilde{\mathbb{P}}^{-y}(\cdot|\xi)$, $\mathbf{m}$ is a Poisson point
process.  Thus the process $\sum_{s\leq
t}1_{\{m_{s}\geq 1\vee e^{N(\xi_{s}+\underline\lambda
s)}\}}$ is a Poisson process with instant intensity $\int_1^\infty
r\nu({\rm d}r) 1_{\{r>e^{N(\xi_t+\underline\lambda t)}\}}{\rm d}t$
under $\widetilde{\mathbb{P}}^{-y}(\cdot|\xi)$ and its domain is
$\{\tau_i: i=1,2,\ldots\}$.  Therefore,
 \begin{eqnarray*}\label{equiv}
&& \sum_{\tau_i<\infty}1_{\left\{m_{\tau_i}\geq
 e^{N(\xi_{\tau_i}+\underline\lambda\tau_i)}\right\}}<\infty\quad
\widetilde{\mathbb{P}}^{-y}(\cdot|\xi)-\mbox{a.s.}\\
&\Longleftrightarrow& \int_0^\infty  {\rm d}t \int_1^\infty r\nu({\rm d}r)
1_{\{r>e^{N(\xi_t+\underline\lambda t)}\}}<\infty\quad
\widetilde{\mathbb{P}}^{-y}(\cdot|\xi)-\mbox{a.s.}
\end{eqnarray*}
Let $a$ be some constant and define the set $C=\left\{\int_0^\infty
 {\rm d}t \int_1^\infty r\nu({\rm d}r) 1_{\{r>e^{N(\xi_t+\underline\lambda
t)}\}}<a\right\}$. Recall that $y+\xi_{t}+\underline\lambda t$ is a
$\mbox{BES}^3(y)$ process under the probability
$\widetilde{\mathbb{P}}^{-y}$.  It is well known that
$\mbox{BES}^3(y)$ is identically distributed to the modulus process
of $B_t+\hat{y}$, where $(B_t, \mathbb Q)$ is a three dimensional
Brownian motion starting at $0$  and $\hat{y}$ is a
point in $\mathbb R^3$ with norm $y$. Denote the modulus process by
$|B_t+\hat{y}|$. We still use $C$ to denote the same set
corresponding to $(B_t, \mathbb Q).$
\begin{eqnarray*}
\lefteqn{\hspace{-2cm}\widetilde{\mathbb{E}}^{-y}\left[1_C\int_0^\infty  {\rm d}t \int_1^\infty
r\nu({\rm d}r) 1_{\{r>e^{N(\xi_t+\underline\lambda t)}\}}\right]}&&\\
&=&\int_0^\infty {\rm d}t\int_1^\infty
r\nu({\rm d}r)\widetilde{\mathbb{E}}^{-y}\left[
1_C1_{\{r>e^{N(\xi_t+\underline\lambda t)}\}}\right]\\
&=&\int_1^\infty r\nu({\rm d}r)\int_0^\infty
\widetilde{\mathbb{E}}^{-y}\left[1_C1_{\{\xi_t+\underline\lambda
t\leq
N^{-1}\log r\}}\right] {\rm d}t \\
&=&\int_1^\infty r\nu({\rm d}r)\int_0^\infty \mathbb
Q\left[1_C1_{\{|B_t+\hat{y}|\leq y+
N^{-1}\log r\}}\right] {\rm d}t \\
&\geq&\mathbb Q\left[1_C\int_1^\infty r\nu({\rm d}r)\int_0^\infty
1_{\{|B_t|\leq N^{-1}\log r\}} {\rm d}t \right].
\end{eqnarray*}
Then under $\mathbb Q$, $|B_t|$ is a $BES^3(0)$ process.  Let
$l_\infty^a$ be the local time of $|B_t|$.
Exercise $(2.5)$ in
\cite{R} tells us $l_\infty^a$ is a $\mbox{BESQ}^2(0)$. Then
$l_\infty^a\stackrel{d}{=}al_\infty^1$ and $\mathbb
Q(l_\infty^1=0)=0.$ For the given set $C$,
\begin{eqnarray*}
\lefteqn{\hspace{-2cm}\mathbb Q\left[1_C\int_1^\infty r\nu({\rm d}r)\int_0^\infty 1_{\{|B_t|\leq
N^{-1}\log r\}} {\rm d}t \right]}&&\\
&&=\mathbb Q\left[1_C\int_1^\infty
r\nu({\rm d}r)\int_0^{N^{-1}\log r}l_\infty^a  {\rm d}a \right]\\
&&=\mathbb Q\left[1_C\int_0^\infty l_\infty^a  {\rm d}a \int_{e^{Na}}^\infty
r\nu({\rm d}r)\right]\\
&&=\mathbb Q\left[1_C\int_0^\infty a  {\rm d}a
\int_0^{a^{-1}l_\infty^a} {\rm d}u \int_{e^{Na}}^\infty
r\nu({\rm d}r)\right]\\
&&=\int_0^\infty a  {\rm d}a \int_{e^{Na}}^\infty r\nu({\rm d}r)\int_0^\infty
\mathbb Q\left[1_C1_{\{l_\infty^a>au\}}\right] {\rm d}u \\
&&\geq\int_0^\infty a {\rm d}a \int_{e^{Na}}^\infty r\nu({\rm d}r)\int_0^\infty
\left[\mathbb Q(C)-\mathbb Q(a^{-1}l_\infty^a>u)\right]^+ {\rm d}u  \\
&&=\int_0^\infty a  {\rm d}a \int_{e^{Na}}^\infty r\nu({\rm d}r)\int_0^\infty
\left[\mathbb Q(C)-\mathbb Q(l_\infty^1<u)\right]^+ {\rm d}u.
\end{eqnarray*}
Note that
\[
\int_0^\infty a  {\rm d}a \int_{e^{Nt}}^\infty r\nu({\rm d}r)=\int_1^\infty
r\nu({\rm d}r)\int_0^{N^{-1}\log r}a {\rm d}a  =\frac{1}{2N^2}\int_1^\infty
r(\log r)^2\nu({\rm d}r).
\]
Therefore, if $\widetilde{\mathbb P}^{-y}(C)>0$, then $\int_1^\infty
r(\log r)^2 \nu({\rm d}r)<\infty$, which means $\int_1^\infty
r(\log r)^2\nu({\rm d}r)=\infty$ implies
$m_{\tau_i}>e^{N(\xi_{\tau_i}+\underline\lambda\tau_i)}$ infinitely
times $\widetilde{\mathbb P}^{-y}$-a.s.  The process
$\mbox{BES}^3(y)$ is transient, so
$\lim_{t\rightarrow\infty}\xi_t+\underline\lambda t=\infty$
$\widetilde{\mathbb P}^{-y}$-a.s.  We reach the conclusion that for
any $N>0,$ there exist an increasing sequence of stochastic time
$\{\tau_i:i=1,2,\ldots\}$ such that
\begin{eqnarray}\label{cri-infsup}
\limsup_{i\rightarrow\infty}m_{\tau_i}
(\xi_{\tau_i}+\underline\lambda\tau_i)e^{-N(\xi_{\tau_i}+\underline\lambda\tau_i)}=\infty.
\end{eqnarray}
Consider the process $V^{\Lambda, -y}_t(\underline\lambda)$ in
\eqref{cri-decom}.  We deduce from \eqref{cri-infsup} that
\begin{eqnarray*}
&&\limsup_{t\rightarrow\infty} V^{\Lambda, -y}_t(\underline\lambda)
\geq\limsup_{i\rightarrow\infty}V^{\Lambda,
-y}_{\tau_i}(\underline\lambda)\\
&\geq&\limsup_{i\rightarrow\infty}\sum_{s\leq \tau_i:\mathbf{m}}
\frac{(y + \xi_s + \underline{\lambda} s)}{y}e^{-\underline
\lambda(\xi_s +\underline\lambda s)}V_{\tau_i-s}^{{\mathbf m}, s,
-(y + \xi_s + \underline{\lambda} s)}(\underline\lambda)\\
&\geq&\limsup_{i\rightarrow\infty}y^{-1}m_{\tau_i}(y + \xi_{\tau_i}
+ \underline{\lambda} \tau_i)e^{-\underline \lambda(\xi_{\tau_i}
+\underline\lambda \tau_i)}\\
&=&\infty.
\end{eqnarray*}
 \qed

\section{Proof of Theorem \ref{main-2}}

(i) Under the given conditions, we know that there exists an $L^1(\mathbb{P})$-limit $W_\infty(\lambda)$ for the martingale $W(\lambda)$. In light of Corollary \ref{support} we have through, a now familiar projection, that
\[
\mathbb{E}(W_\infty(\lambda)|\mathcal{F}^{c_{\lambda}}_{D_x})
= \lim_{t\uparrow\infty}\mathbb{E}(W_\infty(\lambda)|\mathcal{F}^{c_{\lambda}}_{D^t_x})  = e^{-\lambda x}Z^{c_\lambda}_x. \]
This has the
implication that the normalizing sequence discussed in the proof of
Theorem \ref{main-1} (ii) must satisfy $L_\lambda(e^{-\lambda
x})\sim 1$ as $x\uparrow\infty$ and, up to a non-negative constant,
$\Delta(\lambda) = W_{\infty}(\lambda)$.

(ii) In a similar fashion, we note that under the given conditions of the theorem, for fixed $y>0$, $V^{-y}_\infty$ is an $L^1(\mathbb{P})$-limit and hence
\[
\mathbb{E}(yV_\infty^{-y}|\mathcal{F}^{c_{\underline\lambda}}_{D_x}) =(y+x)e^{-{\underline\lambda} x} Z^{c_{\underline\lambda}}_x.
\]
It follows that the normalizing sequence discussed in the proof of Theorem \ref{main-1} (ii) must instead satisfy $L_\lambda(e^{-\underline\lambda x})\sim x$ as $x\uparrow\infty$ and, taking account of (\ref{V=W}), we have, up to a non-negative constant, $\Delta(\underline\lambda) = \partial W_{\infty}(\underline\lambda)$.

\section*{Acknowledgements}
AEK would like to thank Thomas Duquesne for an inspiring discussion.  The research of YXR  is supported in part by NNSF of China
 (Grant No. 10871103 and 10971003) and Specialized Research Fund for the Doctoral Program of Higher Education. The research of AMS is supported by CONACyT grant number 000000000093984.
All four authors are deeply indebted to two referees and an associate editor whose detailed analysis of an earlier version of this paper made for significant improvements.


\begin{thebibliography}{99}

\bibitem{AW} D.G. Aronson and H.F.  Weinberger (1978): Multidimensional nonlinear diffusion arising in population genetics. {\it Adv. in Math.} {\bf 30},  33--76.


\bibitem{BK} J.D. Biggins and A.E. Kyprianou (2004): Measure
change in multitype branching, {\em Adv. Appl. Probab.} {\bf
36}, 544-581.

\bibitem{Bra} M. Bramson (1983): Convergence of solutions of the Kolmogorov equation to travelling waves. {\it Mem. Amer. Math. Soc.} {\bf 44}, no. 285, iv+190 pp.



\bibitem{chauvin}  B. Chauvin (1991) Multiplicative martingales and stopping lines
for branching Brownian motion. \textit{Ann. Probab.} \textbf{30} 1195--1205.


\bibitem{Durrett}  R. Durrett (1996): Probability theory and examples (second edition).{\em Duxbury Press}.

\bibitem{DN} R. Durrett and C.  Neuhauser (1994): Particle systems and reaction-diffusion equations. {\it Ann. Probab.} {\bf 22},  289--333.

\bibitem{Dyn1991} E.B. Dynkin (1991): A probabilistic approach to one class of non-linear differential equations. {\it Probab. Th. Rel. Fields} {\bf 89}, 89--115.

\bibitem{D} E.B. Dynkin (1991):
         Branching particle systems and superprocesses, {\em  Ann. Probab.} {\bf 19(3)}, 1157-1194.

\bibitem{Dyn1993} E.B. Dynkin (1993):  Superprocesses and partial differential equations.
             {\em Ann. Probab.} {\bf 29},1833-1858.


\bibitem{Dyn2001} E.B. Dynkin (2001): Branching exit Markov systems and
superprocesses.
{\it Ann. Probab.} {\bf 29}, 1833--1858.

\bibitem{Dyn2002} E.B. Dynkin (2002): Diffusions, Superdiffusions and Partial Differential Equations.
                AMS, Providencem R.I.

 \bibitem{Dynkin-Kuznetsov} E.B. Dynkin and S.E. Kuznetsov (2004):  $\mathbb{N}$-measures for branching Markov exit systems and their applications to differential equations. {\it Probab. Theory Relat. Fields.} {\bf 130}, 135-150.


\bibitem{EK} J.  Engl\"{a}nder  and A.E. Kyprianou (2004): Local
extinction
versus local exponential growth for spatial
branching processes. {\it  Ann. Probab.} {\bf  32}, 78--99.



\bibitem{EP}  J. Engl\"ander and R.G. Pinsky (1999): On the construction and support properties of measure-valued
                diffusions on $D^d$ with spatially dependent branching.
              {\em  Ann. Probab.} {\bf 27(2)},
              684-730.

              \bibitem{Evans} S.N. Evans (1993):  Two representations of a superprocess. {\it Proc. Royal. Soc. Edin.} {\bf 123A} 959-971.

\bibitem{FM} P.C. Fife and J.B. McLeod, (1977): The approach of solutions of nonlinear diffusion equations to travelling front solutions. {\it  Arch. Ration. Mech. Anal.} {\bf 65},  335--361.


\bibitem{Fish} R.A.   Fisher (1937): The advance of advantageous genes. \textit{%
Ann. Eugenics.} \textbf{7}, 355--369.

              \bibitem{Fitz}  P.J. Fitzsimmons (1988): Construction and regularity of measure-valued Markov branching processes. {\it Israeli J. Math.} {\bf 64}, 337--361.

\bibitem{GHH} Y. Git, J.W. Harris and S.C. Harris (2007): Exponential growth rates in a typed branching diffusion. {\it Ann. Appl. Probab.} {\bf 17}, 609-653.

\bibitem{Grey}  D.R. Grey (1974): Asymptotic behavior of continuous time, continuous state-space branching processes.
              {\em J. App. Probab.} {\bf 11},
              669-677.

\bibitem{HH} R. Hardy and S.C. Harris (2009): A spine approach to branching diffusions with applications to Lp-convergence of martingales. {\it S\'eminaire de Probabilit\'es}, XLII, 281-330.

\bibitem{SHarris1999} S.C. Harris (1999): Travelling waves for the F-K-P-P equation via probabilistic arguments.
              {\em Proc. Roy. Soc. Edin.} {\bf 129(A)},
              503-517.



\bibitem{HR} S.C. Harris and M. Roberts (2009): Measure changes with extinction. {\em Stats. Prob. Letters.} {\bf 79},  1129-1133.


\bibitem{Kam} Y. Kametaka. (1976) On the nonlinear diffusion equation of Kolmogorov-Petrovskii-Piskunov type. {\it Osaka J. Math.} {\bf 13},  11--66.


\bibitem{Kol} A. Kolmogorov, I. Petrovskii and N. Piskounov (1937):
\'{E}tude de l'\'{e}quation de la diffusion avec croissance de la
quantit\'{e} de la mati\`{e}re at son application a un probl\`{e}m
biologique. \textit{Moscow Univ. Bull. Math.} \textbf{1}, 1--25.


\bibitem{Kyprianou2004} A.E. Kyprianou (2004): Travelling wave solution to the
               K-P-P equation: alternatives to Simon
               Harris'probabilistic analysis
               {\it Ann. Inst. H. Poincar\'e}. {\bf 40}, 53-72.
              503-517.


\bibitem{K1} A.E. Kyprianou (2005): Asymptotic radial speed of the support of supercritical branching Brownian motion and
              super-Brownian motion in $\R^d$,
              {\em Markov Proc. Relat. Fields}. {\bf 11},
              145-156.


\bibitem{K} A.E. Kyprianou (2006): {\it Introductory lectures on fluctuations of L\'evy processes with applications.} Springer.


\bibitem{Lp} A.E. Kyprianou and A. Murillo-Salas (2011): Super-Brownian motion: $L^p$-convergence of martingales through the pathwise spine decomposition. {\it Preprint.}

\bibitem{Lau} K.-S. Lau (1985):
On the nonlinear diffusion equation of Kolmogorov, Petrovsky, and Piscounov.
{\it J. Differential Equations} {\bf 59},  44--70.

\bibitem{LG99} J.F. Le Gall (1999): {\it Spatial branching processes, random snakes and partial differential equations.}Lectures in Mathematics ETH Z\"urich. Birkh\"auser Verlag, Basel.

 \bibitem{Rong-Liet.al.} R.-L. Liu, Y.-X. Ren and R. Song (2009): $L \log L$ criterion for a class of superdiffusions.    {\em J. Appl. Probab.} {\bf 46}, 479-496.




\bibitem{Ly}  R. Lyons (1997): A simple path to Biggins' martingale
convergence theorem. In Classical and Modern Branching Processes (K.B.\
Athreya and P. Jagers, eds) 84, 217--222. Springer, New York.

\bibitem{LPP}  R. Lyons, R. Pemantle and Y. Peres (1995): Conceptual proofs of
$L\log L$ criteria for mean behaviour of branching processes. \textit{Ann.
Probab.} \textbf{23}, 1125--1138.

\bibitem{Mal} P. Maillard (2011): The number of absorbed individuals in branching Brownian motion with a barrier. \texttt{arXiv:1004.1426}

\bibitem{McKean}  H.P. McKean (1975): Application of Brownian motion to the
equation of Kolmogorov-Petrovskii-Piskunov.\ \textit{Comm. Pure Appl. Math.}
\textbf{XXIX}, 323--331


\bibitem{neveu}  J. Neveu (1988): Multiplicative martingales for spatial
branching processes. In Seminar on Stochastic Processes 1987, eds E. \c{C}%
inlar, K.L. Chung, R.K. Getoor. Progress in Probability and Statistics, 15,
223--241. Birkha\"{u}ser, Boston.


\bibitem{P} R.G. Pinsky (1995): K-P-P-type asympotics for nonlinear diffusion in a large ball with infinite boundary data and on $\bold R^d$ with infinite initial data outside a large ball. {\it Comm. Partial Differential Equations} {\bf 20},  1369--1393.



  %


\bibitem{RW} Y.-X. Ren and H. Wang (2008): On states of total weighted occupation times of a class of infinitely
divisible superprocesses on a bounded Domain. {\em Potential Anal.} {\bf 28}, 105-137.

\bibitem{R} D. Revuz and M. Yor (1980): Continuous martingales and Brownian motion.
{\em Grundlehren der mathematischen Wissenschaften}, {\bf 293}.
Springer-Verlag, Berlin Heidelberg, New York.


\bibitem{Sheu1997}  Y.C. Sheu (1997):  Lifetime and compactness of range for $\psi$-super-Brownian motion with a general
                branching mechanism.
              {\em Stoch. Proc. Appl.} {\bf 70},
              129-141.

\bibitem{U78}  K. Uchiyama (1978): The behavior of solutions of some fnon-linear diffusion equations for large time. {\em J. Math. Kyoto Univ.} {\bf 18(3)}, 453-508.

\bibitem{Volpert}  A.I. Volpert, V.A. Volpert and V.A. Volpert (1994): Traveling wave solutions of parabolic systems.
{\em Translations of Mathematical Monographs}, {\bf 140}. American
Mathematical Society.


\bibitem{watanabe1968} S. Watanabe (1968): A limit theorem of branching processes and continuous-state branching processes. {\em J. Math. Kyoto Univ.} {\bf 8}, 141--167.
\end{thebibliography}
\end{document}